\documentclass[preprint,12pt]{elsarticle}
\usepackage{hyperref}
\usepackage{color}
\usepackage{amsmath}
\usepackage{amsfonts}
\usepackage{amssymb}
\usepackage{latexsym}
\usepackage{latexsym}
\usepackage{array}
\usepackage{amsthm}
\usepackage{graphics}
\usepackage{comment}
\usepackage{epsf,graphicx,subfigure}
\usepackage{epsfig}
\usepackage{epstopdf}
\usepackage{tikz}

\usetikzlibrary{arrows,positioning,shapes.geometric}

\usepackage[margin=34mm]{geometry}
\usepackage{ulem}
\usepackage{longtable}
\usepackage{rotating}
\usepackage{booktabs}

\newtheorem{Theorem}{Theorem}

\newtheorem{Lemma}{Lemma}
\newtheorem{Remark}{Remark}
\newtheorem{Corollary}{Corollary}
\newtheorem{Example}{Example}
\newtheorem{Assumption}{Assumption}

\newtheorem{Definition}{Definition}

\makeatletter
\def\ps@pprintTitle{%
  \let\@oddhead\@empty
  \let\@evenhead\@empty
  \let\@oddfoot\@empty
  \let\@evenfoot\@oddfoot
}

\makeatother

\begin{document}

\begin{frontmatter}

\title{Convergence Rate Analysis in Limit Theorems for Nonlinear Functionals of the Second Wiener Chaos}

\author[GRLiu]{Gi-Ren Liu\corref{mycorrespondingauthor}}
\cortext[mycorrespondingauthor]{Gi-Ren Liu}
\ead{girenliu@ncku.edu.tw}

\address[GRLiu]{Department of Mathematics, National Cheng Kung University, Tainan, Taiwan
\\ National Center for Theoretical Sciences, National Taiwan University, Taipei, Taiwan}

\begin{abstract}
This paper analyzes the distribution distance between random vectors from the analytic wavelet transform of squared envelopes of Gaussian processes and their large-scale limits. For Gaussian processes with a long-memory parameter $\beta$ below 1/2, the limit combines the second and fourth Wiener chaos. Using a non-Stein approach, we determine the convergence rate in the Kolmogorov metric. When the long-memory parameter $\beta$ exceeds 1/2, the limit is a chi-distributed random process, and the convergence rate in the Wasserstein metric is determined using multidimensional Stein's method. Notable differences in convergence rate upper bounds are observed for long-memory parameters within (1/2,3/4) and (3/4,1).
\end{abstract}

\begin{keyword}
analytic wavelet transform;
central limit theorems;
long-range dependence;
Malliavin calculus;
multidimensional Stein's method;
non-central limit theorems;
rate of convergence.
\MSC[2010] Primary 60G60, 60H05, 62M15; Secondary 35K15.

\end{keyword}

\end{frontmatter}



\section{Introduction}
This paper investigates the distribution of a random process derived from a composition of local and nonlocal transformations applied to a stationary random process. The composition, denoted as $U[s]$ and illustrated in Figure \ref{FlowChart}, involves the complex modulus and
the analytic wavelet transform at the scale $s>0$ formulated via the Hilbert transform \cite{chaudhury2010hilbert,lilly2010analytic,mallat1999wavelet}.
Both the transformation $U[s]$ and its square $U^{2}[s]$, often referred to as the scalogram, have a wide range of applications in signal processing \cite{liu2020diffuse,sejdic2008quantitative} and play a significant role in memory parameter estimation \cite{bruna2015intermittent, clausel2014wavelet,moulines2007spectral}.
For the purpose of establishing a mathematical foundation for the convolutional neural network,
the transformations $\{U[s]\}_{s>0}$ are conceived as neurons, forming the second-order scattering network \cite{guth2021phase,mallat2012group}.
In this network, the output of a neuron $U[s']$, where $s'>0$, serves as input for another $U[s]$.
To streamline numerical computation and theoretical analysis,
an alternative set of neurons $\{U^{2}[s]\}_{s>0}$, utilizing the squared complex modulus as the activation function, is explored in \cite{balestriero2017linear, lostanlen2021one}.
Based on this context, we consider a chi-square random process $Y$, derived from the squared magnitude of
the analytic counterpart of a real-valued stationary Gaussian process $X$, to serve as input for $U[s]$ and its squared form; see Section \ref{sec:envelopeGaussian} for details.

Analyzing the distribution of the process $U^{2}[s]Y:=\{U^{2}[s]Y(t)\mid t\in \mathbb{R}\}$ provides insights that aid in understanding the distribution of outputs from
$U^{2}[s]$ applied to a clean signal in the presence of noise.
However, it is challenging to formulate a clear expression for the distribution of $U^{2}[s]Y$.
Recognizing this limitation, researchers turn their attention to exploring large-scale limits of various functionals of random processes \cite{alodat2020limit,clausel2012large,dobrushin1979non,ivanov2013limit}.
Consequently, evaluating the convergence rate has become a focal point for experts \cite{azmoodeh2022optimal,MR1642664,nourdin2009stein,nourdin2012normal}.
\tikzstyle{line} = [draw, -latex']
\tikzstyle{arrow} = [thick,->,>=stealth]
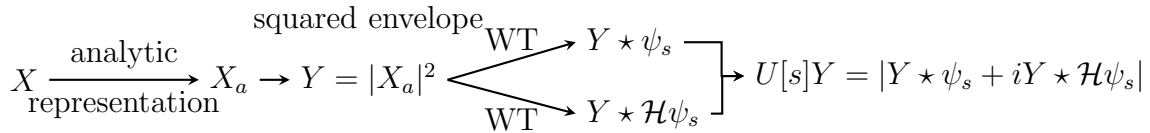
\begin{figure}[h!]  
\centering
   \begin{tikzpicture}[>=latex']
        \tikzset{
        block/.style= {draw, rectangle, align=center,minimum width=1.1cm,minimum height=.10cm,line width=0.05mm},
        rblock/.style={draw, shape=diamond,rounded corners=1.5em,align=center,minimum width=3cm,minimum height=.10cm},
        }
        \node []  (1G) {$X$};
        \node [, right =1.99cm of 1G]  (2F) {$X_{a}$};

        \node [, right =0.40cm of 2F]  (2G) {$Y = |X_{a}|^{2}$};
        \node [, above =0.1 of 2G] (2R) {squared envelope} ;
        \node [, right =1.705cm of 2G, yshift = 0.5cm]  (3F1) {$Y\star \psi_{s}$};
        \node [, right =1.703cm of 2G, yshift = -0.45cm]  (3F2) {$Y\star \mathcal{H}\psi_{s}$};
        \node [, right =0.8cm of 3F1, yshift = -0.45cm]  (3G) {$U[s]Y =|Y\star \psi_{s}+iY\star \mathcal{H}\psi_{s}|$};

          \draw[arrow]   (1G)   --node[above]{analytic} node[below]{representation}(2F); 
          \draw[arrow]   (2F.east)  --(2G.west) node [pos=0.1,above] {}; 
          \draw[arrow]   (2G.east) --(3F1.west) node [pos=0.5,above] {WT};
          \draw[arrow]   (2G.east)  --(3F2.west) node [pos=0.5,below] {WT};
          \draw[arrow]   (3F1.east)  - ++ (0.5cm,0cm) |-(3G.west) node [pos=0.7,above] {};
          \draw[arrow]   (3F2.east)   - ++ (0.15,0cm) |-(3G.west) node [pos=0.8,below] {};

    \end{tikzpicture}
    \caption{The relationship among notations frequently used in this paper: $X$ denotes the Gaussian input process, WT represents the wavelet transform, $\star$ is the convolution operator, $\mathcal{H}$ stands for the Hilbert transform, and $U[s]Y$ or its squared signifies the output.
        }\label{FlowChart}
\end{figure}

The main contributions of this work focus on the large-scale limits of the process $U[s]Y$ and their convergence rates.
In Theorem \ref{thm:nonStein}, we show that when the initial input $X$ is relatively strong long-range dependent, specifically with a long-memory parameter less than 1/2,
the large-scale limit of $U^{2}[s]Y$ has the form of an additive combination of the second and fourth Wiener chaos.
Furthermore, we apply the analytical probability method introduced in the work \cite{MR4003569} by Anh, Vo, Nikolai Leonenko, Andriy Olenko, and Volodymyr Vaskovych  to obtain the convergence rate in the uniform (Kolmogorov) metric.
In Theorem \ref{thm:Stein}, we show that when the initial input $X$ is relatively weak long-range dependent, specifically with a long-memory parameter exceeding 1/2,
the large-scale limit of $U^{2}[s]Y$
is a correlated chi-square distributed random process.
For this case, we apply the work \cite{nourdin2009stein,nourdin2012normal} by Nourdin and Peccati
to calculate the convergence rate
in the Wasserstein metric.

It's noteworthy that, as indicated in the literature \cite{anh1999non,dobrushin1979non,ivanov2013limit,taqqu1979convergence},
some relatively weaker long-range dependent stochastic processes (referred to here as long-memory parameters between 1/2 and 1) may transition into short-range dependent stochastic processes after undergoing local or non-local transformations,
and a distinction exists between the large-scale limits of functionals of short-range dependent random processes and those corresponding to long-range dependent ones.
Therefore, obtaining different large-scale limits under the conditions of Theorems \ref{thm:nonStein} and \ref{thm:Stein}, respectively, is anticipated.
From Theorem \ref{thm:Stein}, we further observe that when the scale variable $s$ tends to infinity,
the upper bound of the Wasserstein distance decreases from $s^{-1/2}L^{2}(s)$ to $s^{0}L^{2}(s)$ following the function $s^{-1/2+2(3/4-\beta)}L^{2}(s)$
when the strength of long-range dependence of the initial input $X$ increases (i.e.,  the long-memory parameter $\beta$ of $X$ decreases from 3/4 to 1/2),
where $L$ is a slowly varying function.
To the best of our knowledge,
this relationship between the rate of convergence and the strength of long-range dependence in initial inputs has hardly been observed in the existing literature concerning
the central limit theorem arising from convolutional functionals of random processes.
Our findings resonate with those in the literature \cite{bierme2011central,breton2008error}
about the total variation distance between the normalized sum of the Hermite power of increments of a fractional Brownian motion and its limit,
as well as how the rate of convergence depends on the Hurst index of the fractional Brownian motion.

In contrast to our previous works \cite{liu2022asymptotic,liu2023central,liu2021scattering},
we consider a transformation with more practical significance in this work.
Specifically, we replace the conventional real-valued wavelet transform with the analytic wavelet transform formulated via the Hilbert transform.
The resulting output is a complex-valued field indexed by time and scale variables.
At each fixed scale, the magnitude of this complex-valued field serves as the envelope of the wavelet coefficient process, providing enhanced stability against time shifts in the underlying time-evolution process. This transform has been extensively applied to extract rich information from physiological signals \cite{liu2020diffuse} and audio recordings \cite{anden2011multiscale,balestriero2017linear} for intelligent classification.
Because the real and imaginary parts of the output of the analytic wavelet transform are correlated
and we need to deal with the sum of their individual square,
existing results associated with the real-valued wavelet transform cannot be applied to the complex-valued case.
Furthermore, we employ both Stein and non-Stein approaches to calculate the upper bounds for the rates of convergence in the Kolmogorov and Wasserstein metrics of the outputs of $U[s]$ and its square to their large-scale limits.
It is worth mentioning that
we verified the positive definiteness condition of the covariance matrix of the normal random vector employed in the multivariate normal approximations. Different from the smooth Wasserstein metric considered in \cite{liu2022asymptotic,liu2023gaussian}, this crucial verification enables us to use the Wasserstein metric
to gauge the distance between the distributions of the outputs of $U[s]$ and their large-scale limits, as shown in Theorem \ref{thm:Stein}.
Moreover, we provide an in-depth analysis of how these upper bounds vary with the long-memory parameter of the initial input $X$.

The article is organized as follows. In Section \ref{sec:preliminary}, we recall some definitions and formulae related to
the analytic wavelet transform and
the envelope of stationary Gaussian processes.
In Section \ref{sec:mainresult}, we state our main
results, including Theorems \ref{thm:nonStein}
and \ref{thm:Stein}, along with key lemmas. Some conclusions are drawn in Section \ref{sec:Conclusion}.
The proofs of our main results are given
in Section \ref{sec:proof}.
Additionally,
the appendix includes lemmas concerning slowly varying functions and anti-concentration inequalities.

\section{Preliminaries}\label{sec:preliminary}
\subsection{Analytic wavelet transform}\label{sec:preliminary:wavelet}

In the field of signal segment classification, considering a real-valued process $\{Y(t)\}_{t\in \mathbb{R}}$ and a set of real-valued testing functions $\{\psi_{s}\}_{s\in \Lambda}$, where $\Lambda$ is an index set, for assessing the similarity between the pattern of the process $Y$ during a time neighborhood of $t$ and $\psi_{s}$ for each $s\in \Lambda$, we consider the convolution
\begin{align}\label{def:convolution}
Y\star \psi_{s}(t) = \int_{\mathbb{R}}Y(t-u)\psi_{s}(u)du,\ t\in \mathbb{R}.
\end{align}
In the field of wavelet theory, $\Lambda = (0,\infty)$ and
\begin{align*}
\psi_{s}(u) = \frac{1}{s}\psi\left(\frac{u}{s}\right),
\end{align*}
where $\psi$ is an integrable function defined on $\mathbb{R}$ satisfying $\int_{\mathbb{R}}\psi(u)du=0$.
For small values of $s$, $\psi_{s}$ quantifies the small-scale fluctuations in $Y$ through (\ref{def:convolution}).
Subsequently, we refer to $s$ as the scale variable.
To reduce the impact of small time shifting on the value of convolution (\ref{def:convolution}), one prefers to use the envelope of
$Y\star \psi_{s}$, which can be obtained from the magnitude of the following complex wavelet transform
\begin{align}\label{def:AWT}
\mathcal{W}_{s}Y(t)=Y\star \psi_{s}(t)+i Y\star \mathcal{H}\psi_{s}(t),
\end{align}
where $\mathcal{H}$ represents the Hilbert transform defined by
\begin{align*}
\mathcal{H}\psi_{s}(t) = \frac{1}{\pi}\int_{-\infty}^{\infty}
\frac{\psi_{s}(u)}{t-u}du.
\end{align*}
We denote the squared envelope of precess $Y\star \psi_{s}$ as $U^{2}[s]Y$, i.e.,
\begin{align}\label{def:;scalogram}
U^{2}[s]Y(t):=\left(Y\star \psi_{s}(t)\right)^{2}+\left( Y\star \mathcal{H}\psi_{s}(t)\right)^{2},
\end{align}
which is known as the scalogram.
We denote the Fourier transform of $\psi$ by $\widehat{\psi}$, which is defined as
\begin{align*}
\widehat{\psi}(\lambda) = \int_{\mathbb{R}}e^{-iu\lambda}\psi(u)du,\ \lambda\in \mathbb{R}.
\end{align*}
As stated in \cite{chaudhury2010hilbert}, $\widehat{\psi_{s}}(\lambda)=\widehat{\psi}(s\lambda)$ and $\widehat{\mathcal{H}\psi_{s}}(\lambda)=-i\textup{sgn}(\lambda)\widehat{\psi}(s\lambda)$, where
$\textup{sgn}(\lambda)=\lambda/|\lambda|$ for $\lambda\neq0$ and $\textup{sgn}(0)=0$.
Because $\psi_{s}+i \mathcal{H}\psi_{s}$ is analytic, meaning its Fourier transform is null for negative frequencies, (\ref{def:AWT}) is referred to as the analytic wavelet transform \cite{lilly2010analytic,mallat1999wavelet}.
In the following content, we only require $\psi$ and $\widehat{\psi}$ to belong to $L^{1}$.
Commonly used wavelets such as the Morlet wavelet, the Daubechies wavelet, the Mexican Hat wavelet, and the Morse wavelet are suitable for our purposes, with the exception of the Haar wavelet.

\subsection{Envelope of stationary Gaussian processes}\label{sec:envelopeGaussian}

Consider a real-valued, mean-square continuous, and stationary Gaussian random process $\{X(t)\}_{t\in \mathbb{R}}$ with a mean of zero and covariance function denoted as $C_{X}$.
It is well known
that
there exists a unique nonnegative measure $F_{X}:\mathcal{B}(\mathbb{R})\rightarrow [0,\infty)$ such that $F_{X}(\Delta) = F_{X}(-\Delta)$ for any $\Delta\in\mathcal{B}(\mathbb{R})$
and
\begin{equation*}
C_{X}(t) = \int_{\mathbb{R}}e^{i\lambda t}F_{X}(d\lambda),\ t\in \mathbb{R}.
\end{equation*}
The measure $F_{X}$ is referred to as the spectral measure.
\begin{Definition}[\cite{bingham1989regular}]\label{def:slowlyvarying}
A measurable function $L :(0,\infty) \rightarrow(0,\infty)$ is said to be slowly varying
at infinity if for all $c> 0$,
\begin{align*}
\underset{s\rightarrow\infty}{\lim}\frac{L(cs)}{L(s)} = 1.
\end{align*}
A measurable function $g: (0,\infty)\rightarrow(0,\infty)$ is said to be regularly
varying at infinity of index $\tau\in \mathbb{R}$ if there exists $\tau$ such that for all $c > 0$,
\begin{align*}
\underset{s\rightarrow\infty}{\lim}\frac{g(cs)}{g(s)} = c^{\tau}.
\end{align*}
\end{Definition}
In accordance with the literature \cite[Assumption 2]{MR4003569}, we make the following assumption on the spectral measure.
\begin{Assumption}\label{Assumption:spectral}
The spectral measure $F_{X}$ is absolutely continuous with respect to the Lebesgue measure with a density function
$f_{X}\in L^{1}(\mathbb{R})$. It has the form
\begin{equation}\label{condition_f}
f_{X}(\lambda) = L(|\lambda|^{-1})|\lambda|^{\beta-1},\ \lambda\neq0,
\end{equation}
where $\beta\in(0,1)$ is the long-memory parameter, $L:(0,\infty)\rightarrow(0,\infty)$ is a locally bounded function that is slowly varying at infinity, and there exists an absolute constant $C$ such that
\begin{align}\label{errorterm_slowlyvarying}
\left|\frac{L(\eta s)}{L(s)}-1\right|\leq Cg_{\tau}(s)h_{\tau}(\eta),\ \eta\geq 1.
\end{align}
In (\ref{errorterm_slowlyvarying}), $g_{\tau}(\cdot)$ is a regular varying function of index $\tau\leq 0$ with $g_{\tau}(s)\rightarrow 0$ as $s\rightarrow\infty$
and $h_{\tau}: [1,\infty)\rightarrow[0,\infty)$ is defined as follows
\begin{align}\label{def:h_tau}
h_{\tau}(\eta) = \left\{\begin{array}{ll}\ln(\eta) &\ \textup{if}\ \tau=0,\\
\tau^{-1}(\eta^{\tau}-1) &\ \textup{if}\ \tau<0.\end{array}\right.
\end{align}
\end{Assumption}
For the case where $L(s) = \ln(s)$,  we have $\tau=0$ and $g_{\tau}(s) = |\ln(s)|^{-1}$ because
$|L(\eta s)/L(s)-1| = \ln(\eta)/|\ln(s)|$.
\begin{Example}\label{label:example:Linnik}
If $C_{X}$ is the generalized Linnik
covariance function \cite{lim2010analytic}, defined as
\begin{align*}
C_{X}(t)=(1+|t|^{\sigma})^{-\chi},\ \sigma\in(0,2],\ \chi>0,
\end{align*}
then we define $\widetilde{L}(|t|) =(1+|t|^{-\sigma})^{-\chi}$, which is slowly varying at infinity.
Because $C_{X}(t) =\widetilde{L}(|t|) |t|^{-\sigma \chi}$, by the Tauberian theorem \cite[Theorem 4]{leonenko2013tauberian}, if  $0<\sigma \chi<1$,
\begin{align*}
f_{X}(\lambda)\sim \widetilde{L}\left(|\lambda|^{-1}\right)|\lambda|^{\sigma \chi-1}\left[2\Gamma(\sigma \chi)\cos\left(\frac{\sigma \chi\pi}{2}\right)\right]^{-1}
\end{align*}
as $|\lambda|\rightarrow0$. Hence, Condition (\ref{condition_f}) holds with $\beta = \sigma\chi$.
For the function $\widetilde{L}$ in this example, by the mean value theorem,
\begin{align}\label{residual_Bessel}
\left|\frac{\widetilde{L}(\eta s)}{\widetilde{L}(s)}-1\right|\leq \sigma\chi h_{-\sigma}(\eta) g_{-\sigma}(s)
\end{align}
for any $s>0$ and $\eta\geq1$, where
\begin{align*}
g_{-\sigma}(s) = \frac{(1+s^{\sigma})^{\chi}}{s^{\sigma(\chi+1)}},\ s>0.
\end{align*}
\end{Example}

\begin{Example}\label{label:example:Dagum}
If $C_{X}$ is the Dagum covariance function \cite{faouzi2022deep},
defined as
\begin{align}\label{Dagum}
C_{X}(t)=1-(1+|t|^{-\sigma})^{-\chi},\ \sigma\in(0,2],\ \chi>0,
\end{align}
then we define $\widetilde{L}(|t|)=\chi^{-1}|t|^{\sigma}C_{X}(t)$.
By the L'Hospital's rule, $\widetilde{L}(t)\rightarrow 1$ as $|t|\rightarrow\infty$.
Hence, $\widetilde{L}$ is a slowly varying function.
It also implies that the random process $X$
exhibits the long-range dependence when $\sigma\in(0,1)$.
Under the condition $\sigma\in(0,1)$ and $\sigma\chi\in (0,2]$,
\cite[Theorem 4]{faouzi2022deep} proves that Condition (\ref{condition_f}) holds with $\beta = \sigma$.
For the function $\widetilde{L}$ in this example, by the mean value theorem, there exists a constant $C$ depending on
$\sigma$ and $\chi$ such that
\begin{align}\label{residual_Dagum}
\left|\frac{\widetilde{L}(\eta s)}{\widetilde{L}(s)}-1\right|\leq Ch_{-\sigma}(\eta) s^{-\sigma}
\end{align}
for any $s\geq1$ and $\eta\geq1$.
For the case $s\in(0,1]$, the graph of the function $\eta \mapsto s^{\sigma}|\widetilde{L}(\eta s)/\widetilde{L}(s)-1|$
is illustrated in Figure \ref{fig:Dagum}, considering $s=10^{-7}, 10^{-6},..., 10^{-1},1.$
In Figures \ref{fig:Dagum}(a)-(c), we consider three different sets of parameters for
the Dagum covariance function, including $(\sigma,\chi)=(0.25,1),(0.5,1), (0.75,1)$.
These three figures demonstrate that for each $\eta\geq1$, $s^{\sigma}\left|\widetilde{L}(\eta s)/\widetilde{L}(s)-1\right|$ monotonically increases with $s$.
This observation implies that
\begin{align*}
s^{\sigma}\left|\frac{\widetilde{L}(\eta s)}{\widetilde{L}(s)}-1\right|\leq Ch_{-\sigma}(\eta)
\end{align*}
for any $s\in(0,1]$ and $\eta\geq1$.
\begin{figure}[htb]
\centering
\subfigure[][$\sigma = 0.25$]
{\includegraphics[scale=0.53]{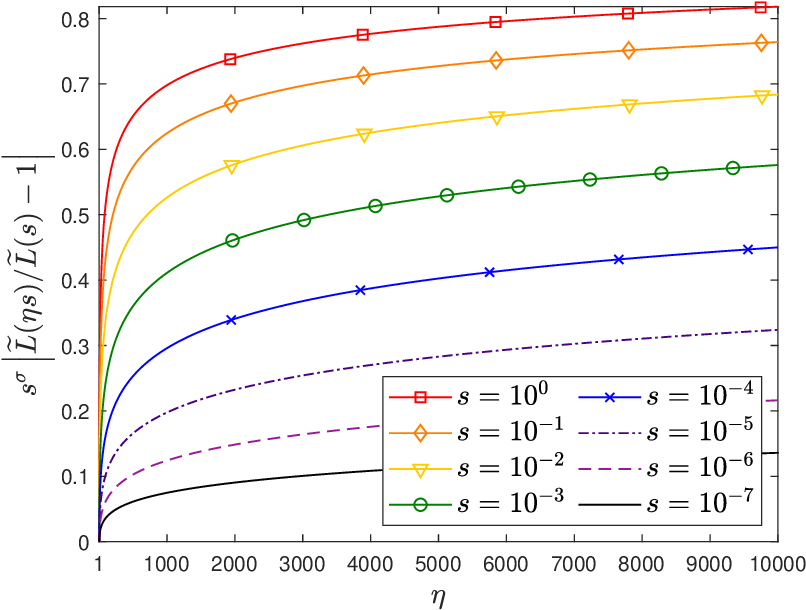}}
\subfigure[][$\sigma=0.50$]
{\includegraphics[scale=0.53]{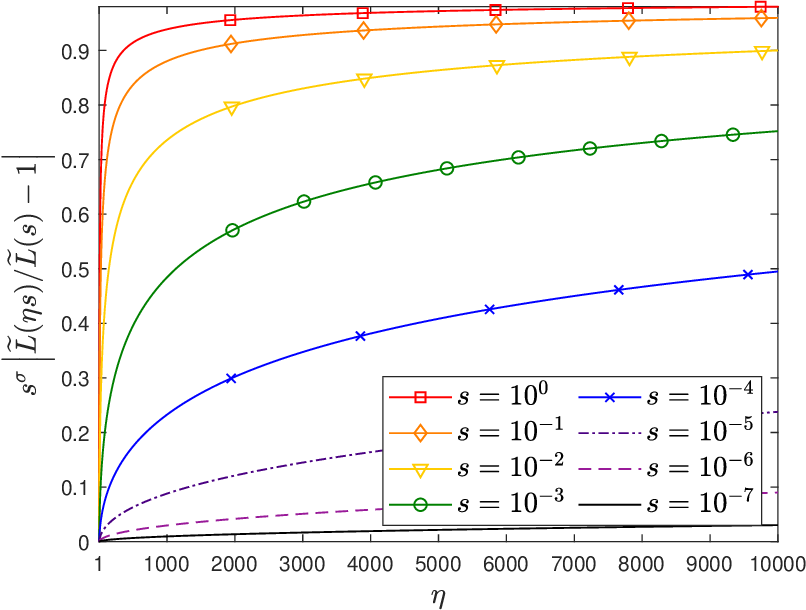}}
\\
\subfigure[][$\sigma=0.75$]
{\includegraphics[scale=0.53]{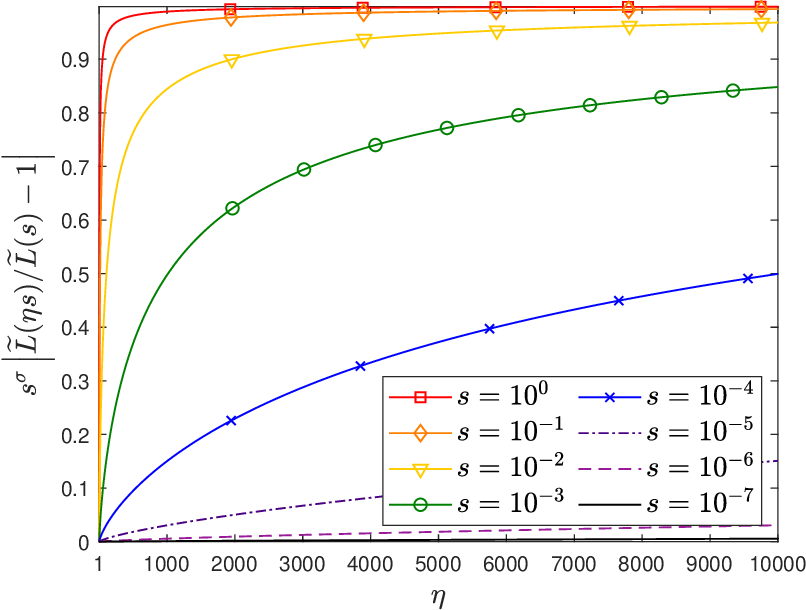}}
\subfigure[][the graph of $h_{\tau}$]
{\includegraphics[scale=0.53]{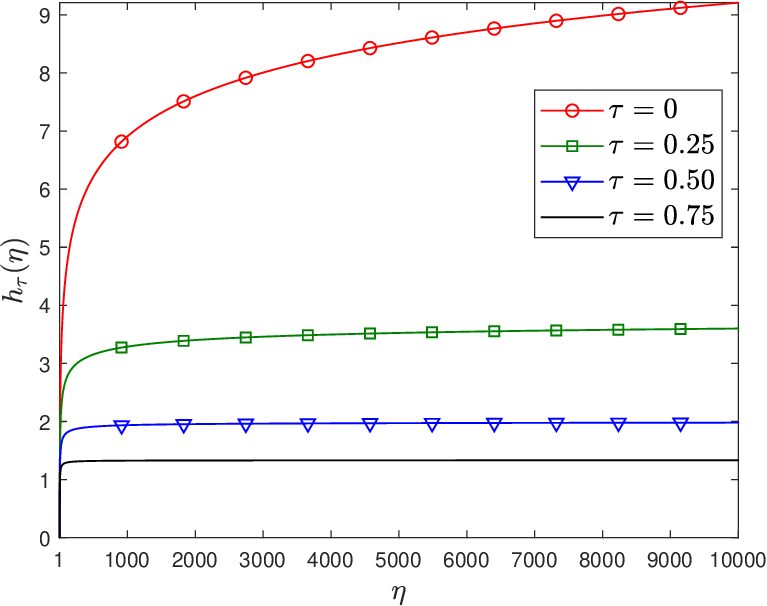}}
\caption{Effects of the variable $s$ and the parameter $\sigma$ on the behavior of the function $\eta\mapsto s^{\sigma}|\widetilde{L}(\eta s)/\widetilde{L}(s)-1|$, where $\widetilde{L}(|t|)=\chi^{-1}|t|^{\sigma}C_{X}(t)$, $C_{X}$ is the Dagum covariance function with parameters $\sigma$ and $\chi$ defined in (\ref{Dagum}). Here, $\chi$ is fixed at 1.}
\label{fig:Dagum}
\end{figure}
\end{Example}

Under Assumption \ref{Assumption:spectral}, and according to the Karhunen theorem,
the Gaussian process $X$ can be represented as follows
\begin{align*}
X(t)=\int_{\mathbb{R}}
e^{i\lambda t}\sqrt{f_{X}(\lambda)}W(d\lambda),
\end{align*}
where $W$ is a complex-valued Gaussian white noise random measure on $\mathbb{R}$
satisfying
\begin{align}\label{ortho}
W(\Delta_{1})=\overline{W(-\Delta_{1})}\ \ \textup{and}\ \ \mathbb{E}\left[W(\Delta_{1})
\overline{W(\Delta_{2})}\right]=\textup{Leb}(\Delta_{1}\cap\Delta_{2})
\end{align}
for any intervals
 $\Delta_{1},\Delta_{2}$, where Leb is the Lebesgue measure on $\mathbb{R}$ \cite{major1981lecture}.

\begin{figure}
\centering
  \includegraphics[scale=0.375]{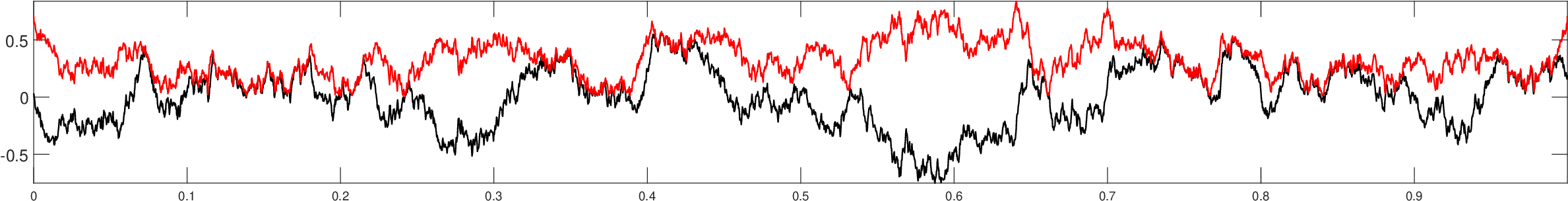}
  \caption{Illustration: a sample path of the Gaussian process $X$ (in black) and its envelope $|X_{a}|$ (in red).}\label{fig:envelope}
\end{figure}

The analytic counterpart of $X$ is a complex-valued process defined by
\begin{align}\notag
X_{a}(t)=& X(t)+i\mathcal{H}X(t)
\\\notag=&
2\int_{0}^{\infty}
e^{i\lambda t}\sqrt{f_{X}(\lambda)}W(d\lambda),\ t\in \mathbb{R},
\end{align}
where the negative frequency content of $X$ is suppressed.
As described in \cite{gabor1946theory,mallat1999wavelet}, from the process $X_{a}$,
we can efficiently extract both the amplitude
and phase components of $X$ through the relationship
$$
X(t) = \textup{Re}(|X_{a}(t)|e^{i\arg(X_{a}(t))}).
$$
In comparison with the phase component $\{\arg(X_{a}(t))\mid t\in \mathbb{R}\}$,
the process $\{|X_{a}(t)|\mid t\in \mathbb{R}\}$, known as the envelope of $X$, is more stable against time shifts.
Figure \ref{fig:envelope} illustrates a sample path of $X$ and its envelope $|X_{a}|$.
This stability can be advantageous in some applications \cite{guth2021phase}.

From another perspective, when convolving a real-valued random process
with an analytic wavelet function, it will generate the analytic counterpart
of the wavelet coefficient process.
This establishes a connection where the envelope process $|X_{a}|$ mentioned earlier can be viewed as the complex modulus
of the wavelet coefficient process of a certain time-dependent process.
In other words, the process $|X_{a}|$ (resp. $|X_{a}|^{2}$) can be interpreted as the output of another neuron $U[\cdot]$ (resp. $U^{2}[\cdot]$) \cite{mallat2012group}.
Based on this context, we consider the scenario
\begin{align}\label{def:envelopeX}
Y(t) = |X_{a}(t)|^{2},\ t\in \mathbb{R},
\end{align}
for the process $Y$ considered in Section \ref{sec:preliminary:wavelet}.
For an analysis of level crossings by the process $Y$, readers can refer to the literature \cite{adler1978envelope}.
To facilitate the subsequent time-frequency analysis, we establish a Wiener integral representation for $Y$ in the following lemma.

\begin{Lemma}\label{lemma:chisquare}
Under the assumption about the existence of the spectral density function $f_{X}$ in Assumption \ref{Assumption:spectral}, $\{Y(t)/\|f_{X}\|_{1}\}_{t\in \mathbb{R}}$ defined in (\ref{def:envelopeX})
is a chi-square distributed random process with two degrees of freedom, and
\begin{align}\label{cov:Y_revised}
\textup{Cov}(Y(t),Y(t')) =16\left|\int_{0}^{\infty} e^{i\lambda(t-t')}f_{X}(\lambda)d\lambda\right|^{2},\ t,t'\in \mathbb{R}.
\end{align}
Moreover, $Y$ has the Wiener integral representation
\begin{align}\label{Ychaosdep}
Y(t) = 2\|f_{X}\|_{1}+2\int_{\mathbb{R}^{2}}^{'}
e^{it(\lambda_{1}+\lambda_{2})}\sqrt{f_{X}(\lambda_{1})f_{X}(\lambda_{2})}1_{\{\lambda_{1}\lambda_{2}<0\}}W(d\lambda_{1})W(d\lambda_{2}),
\end{align}
where $W$ is a complex-valued Gaussian white noise random measure satisfying (\ref{ortho})
and $\int^{'}$ means that the integral excludes the hyperplanes
$\lambda_{k}=\mp \lambda_{k^{'}}$ for $k, k^{'}\in\{1,2\}$ and  $k\neq k^{'}$.

\end{Lemma}

The proof of Lemma \ref{lemma:chisquare} is provided in Section \ref{sec:proof:lemma:chisquare}.
Similar to the Rosenblatt process \cite{dobrushin1979non} and the chi-square random field proposed in
\cite[Example 2.3]{anh1999non},
the process $Y$ is expressed as a double iterated integral with respect to the complex-valued Gaussian white noise random measure $W$.
For the wavelet analysis of the Rosenblatt process, please refer to \cite{bardet2010wavelet,clausel2014wavelet}.
Different from the real-valued wavelets employed in \cite{bardet2010wavelet,clausel2014wavelet},
the wavelets considered in this work are complex-valued and analytic. The reason behind adopting analytic wavelets for signal or random process analysis is their stability in magnitude during local input translations.
This stability proves advantageous in signal classification, as signal classes typically remain invariant to local translations \cite{guth2021phase,mallat2012group}.
The interplay between the real and imaginary components manifests in the spectral representation of the squared magnitude of the analytic wavelet transform of $Y$, as presented in Lemma \ref{lemma:chasodep:UY} below.

\begin{Remark}\label{remark:difference_Y}
In \cite{anh1999non,anh2002renormalization,anh2013macroscaling}, the authors consider
the chi-square random processes
\begin{align}
\chi^{2}_{d}(t) = \frac{1}{2}\left(\zeta_{1}^{2}(t)+\zeta_{2}^{2}(t)+\cdots+\zeta_{d}^{2}(t)\right),\ d\in \mathbb{N},
\end{align}
where $\zeta_{1},\ldots, \zeta_{d}$ are independent copies of a stationary Gaussian process $\zeta$ with
$\mathbb{E}\left[\zeta(t)\right]=0$ and $\mathbb{E}\left[\zeta^{2}(t)\right]=1.$
Upon normalization, the processes $\chi^{2}_{2}$ and  $Y$ follow the same marginal distribution.
However, in the definition of $Y$, i.e., $Y(t) = (X(t))^{2}+(\mathcal{H}X(t))^{2}$, the proof of Lemma 1 shows that
$X(t)$ and $\mathcal{H}X(t')$ are dependent if $t\neq t'$.
This discrepancy leads to differences between the covariance function of $Y$ in (\ref{cov:Y_revised}) and that of $\chi^{2}_{2}$:
$$
\textup{Cov}\left(\chi^{2}_{2}(t),\chi^{2}_{2}(t')\right) = \left(\textup{Cov}\left(\zeta(t),\zeta(t')\right)\right)^{2}.
$$
Additionally, the spectral representation of $Y$ in (\ref{Ychaosdep}) differs from that of $\chi^{2}_{2}$, as shown below:
\begin{align}\notag
\chi^{2}_{2}(t) = & \frac{1}{2}\int_{\mathbb{R}^{2}}^{'}e^{it(\lambda_{1}+\lambda_{2})}f_{\zeta}^{1/2}(\lambda_{1})f_{\zeta}^{1/2}(\lambda_{2})W_{1}(d\lambda_{1})
W_{1}(d\lambda_{2})
\\\notag&+\frac{1}{2}\int_{\mathbb{R}^{2}}^{'}e^{it(\lambda_{1}+\lambda_{2})}f_{\zeta}^{1/2}(\lambda_{1})f_{\zeta}^{1/2}(\lambda_{2})W_{2}(d\lambda_{1})
W_{2}(d\lambda_{2})+1,
\end{align}
where $f_{\zeta}$ is the spectral density of $\zeta$ if it exists, and $W_{1}$ and $W_{2}$ are independent Gaussian white noise random measures.
\end{Remark}

\begin{Lemma}\label{lemma:chasodep:UY}
Consider the random process $Y$ defined in (\ref{def:envelopeX}).
Assuming that $\psi$ and its Fourier transform $\widehat{\psi}$ belong to $L^{1}$, and considering Assumption \ref{Assumption:spectral},
for each $s>0$, the random process $\{U^{2}[s]Y(t)\mid t\in \mathbb{R}\}$ has the Wiener-It$\hat{o}$ decomposition
\begin{align}\notag
U^{2}[s]Y(t) =& \mathbb{E}\left[U^{2}[s]Y(t)\right]+
\int_{\mathbb{R}^{2}}^{'}
e^{it(\lambda_{1}+\lambda_{2})}K_{2}(\lambda_{1},\lambda_{2})
W(d\lambda_{1})W(d\lambda_{2})
\\\notag&+
\int_{\mathbb{R}^{4}}^{'}
e^{it(\lambda_{1}+\lambda_{2}+\lambda_{3}+\lambda_{4})}
K_{4}(\lambda_{1},\lambda_{2},\lambda_{3},\lambda_{4})
\overset{4}{\underset{\ell=1}{\prod}}W(d\lambda_{\ell}),
\end{align}
where
\begin{align*}
\mathbb{E}\left[U^{2}[s]Y(t)\right] = 2^{4}\int_{\mathbb{R}^{2}}|\widehat{\psi}(s(\lambda_{1}+\lambda_{2}))|^{2}
f_{X}(\lambda_{1})
f_{X}(\lambda_{2})1_{\{\lambda_{1}\lambda_{2}<0\}}d\lambda_{1}d\lambda_{2},
\end{align*}
\begin{align*}
K_{2}(\lambda_{1},\lambda_{2}) =& 8\sqrt{f_{X}(\lambda_{1})
f_{X}(\lambda_{2})}1_{\{\lambda_{1}\lambda_{2}<0\}}
\\&\times
\int_{\mathbb{R}}f_{X}(\xi)
\widehat{\psi}(s(\lambda_{1}+\xi))\widehat{\psi}(s(\lambda_{2}-\xi)) 1_{\{(\lambda_{1}+\xi)(\lambda_{2}-\xi)<0\}}d\xi,
\end{align*}
and
\begin{align*}
K_{4}(\lambda_{1},\lambda_{2},\lambda_{3},\lambda_{4}) =&
8\left[\overset{4}{\underset{\ell=1}{\prod}}\sqrt{f_{X}(\lambda_{\ell})}\right]\widehat{\psi}(s(\lambda_{1}+\lambda_{2}))\widehat{\psi}(s(\lambda_{3}+\lambda_{4}))
\\&\times1_{\{\lambda_{1}\lambda_{2}<0\}}1_{\{\lambda_{3}\lambda_{4}<0\}}
1_{\{(\lambda_{1}+\lambda_{2})(\lambda_{3}+\lambda_{4})<0\}}.
\end{align*}
\end{Lemma}

\begin{Remark}
Let $\{L_{k}\}_{k=0,1,2,...}$ denote the Laguerre polynomials, which form an orthogonal basis for the Hilbert space $L^{2}((0,\infty),e^{-r}dr)$. For example, $L_{0}(r)=1$, $L_{1}(r)=1-r$, and $L_{2}(r) = 2^{-1}r^{2}-2r+1$.
The first and second Laguerre chaos $L_{1}(\chi^{2}_{2})$ and $L_{2}(\chi^{2}_{2})$ derived from the chi-square random process
in Remark \ref{remark:difference_Y} can be expressed as linear combinations of the second and fourth Wiener chaos \cite[p. 111]{anh1999non} as follows:
\begin{align}\label{L_H_relation1}
L_{1}(\chi^{2}_{2}(t)) = -\frac{1}{2}\left(H_{2}(\zeta_{1}(t))+H_{2}(\zeta_{2}(t))\right)
\end{align}
and
\begin{align}\label{L_H_relation2}
L_{2}(\chi^{2}_{2}(t)) = \frac{1}{8}\left[2H_{2}(\zeta_{1}(t))H_{2}(\zeta_{2}(t))
-H_{4}(\zeta_{1}(t))-H_{4}(\zeta_{2}(t))\right],
\end{align}
where $H_{k}$, $k=0,1,2,...$, represent the Hermite polynomials.
Several limit theorems, including \cite[Theorem 3.3]{anh1999non},
\cite[Theorem 2.4]{anh2002renormalization} and \cite[Section 5]{anh2013macroscaling}, about additive functionals of the Laguerre chaos have been established based on (\ref{L_H_relation1}), (\ref{L_H_relation2}), and various filtering methods.

It is worth to mentioning that the transformation $U^{2}[s]$ is not an additive functional.
As shown in Figure \ref{FlowChart}, the transformation $U^{2}[s]$ of $Y$ comprises two parallel convolution operators, namely $Y\star\psi_{s}$ and $Y\star \mathcal{H}\psi_{s}$, and the sum of their squares.
The order of the action of local and non-local operators on chi-square random processes is exactly opposite to that in the literature \cite{anh1999non,anh2002renormalization,anh2013macroscaling}.
The discrepancy makes the spectral representation of $U^{2}[s]Y$ in Lemma \ref{lemma:chasodep:UY} difficult to express as Laguerre chaos.
Consequently, the main results presented below cannot be encompassed by existing limit theorems.
\end{Remark}

\section{Main results}\label{sec:mainresult}
In this section, we present our primary findings concerning the large-scale limits of the process $U[s]Y$ and their convergence rates.
We'll observe that the long-memory parameter of the initial input $X$ at $1/2$ serves as a threshold.
When the long-memory parameter falls below or exceeds this threshold, the normalized $U[s]Y$ converges to distinct limits.
Moreover, the convergence rate might vary based on the long-memory parameter, particularly within the interval $(1/2, 3/4)$.

In the subsequent discussion, when comparing two functions, $f(s)$ and $g(s)$, the notation $f(s) \lesssim g(s)$ implies the existence of a constant $C$ that remains unaffected  regardless of changes in the variable $s$, satisfying the condition $f(s) \leq Cg(s)$. The value of constant $C$ may vary between different equations.

\subsection{Non-central limit theorems}

We begin by considering a random process $Z_{s}$ obtained from the squared magnitude of the analytic wavelet transform of  $Y$:
\begin{align}\label{def:Zs}
Z_{s}(t)= s^{2\beta}L^{-2}(s)U^{2}[s]Y(st),
\end{align}
where the time units of $Z_{s}$ are scaled by a factor of $s$ compared to the time units of $Y$.
\begin{Definition}\label{def:Kol:distance}
Consider two random variables $V_{1}$ and $V_{2}$. The Kolmogorov distance between the distributions of $V_{1}$ and $V_{2}$, denoted as $d_{\textup{Kol}}(V_{1},V_{2})$, is defined as follows
\begin{align*}
d_{\textup{Kol}}(V_{1},V_{2}) = \underset{z\in \mathbb{R}}{\sup}\left|\mathbb{P}(V_{1}\leq z)-\mathbb{P}(V_{2}\leq z)\right|.
\end{align*}
\end{Definition}
Theorem \ref{thm:nonStein} below concerns the convergence rate
of the Kolmogorov distance between $Z_{s}$ and its limit as $s\rightarrow\infty$.

\begin{Theorem}\label{thm:nonStein}
Consider the random process $Y$ defined in (\ref{def:envelopeX}). Under the assumption that the wavelet $\psi$ and its Fourier transform $\widehat{\psi}$ belong to $L^{1}$, and considering Assumption \ref{Assumption:spectral} with $0<2\beta<1$ and $-\tau<(1-2\beta)/4$,
for any $n\in \mathbb{N}$ and $c_{1},...,c_{n}, t_{1},...,t_{n}\in \mathbb{R}$, the following result holds:
\begin{align}\label{decay_in_g(s)}
d_{\textup{Kol}}\left(\overset{n}{\underset{k=1}{\sum}}c_{k}Z_{s}\left(t_{k}\right),\overset{n}{\underset{k=1}{\sum}}c_{k}Z(t_{k})\right)\lesssim  \left(g_{\tau}(s)\right)^{2/9},
\end{align}
where $Z(t) =Z^{(0)}+ Z^{(1)}(t)+Z^{(2)}(t)$ with
\begin{align}\label{def:Z0}
Z^{(0)}=16\int_{\mathbb{R}^{2}}|\widehat{\psi}(\lambda_{1}+\lambda_{2})|^{2}
|\lambda_{1}\lambda_{2}|^{\beta-1}1_{\{\lambda_{1}\lambda_{2}<0\}}d\lambda_{1}d\lambda_{2},
\end{align}
\begin{align}\notag
Z^{(1)}(t)=&8\int_{\mathbb{R}^{2}}^{'}
e^{it(\lambda_{1}+\lambda_{2})}\left|\lambda_{1}\lambda_{2}\right|^{\frac{\beta-1}{2}}
1_{\{\lambda_{1}\lambda_{2}<0\}}
\\\label{def:Zlimit1}&\times
\left[\int_{\mathbb{R}}
\frac{\widehat{\psi}(\lambda_{1}+\xi)\widehat{\psi}(\lambda_{2}-\xi)}{\left|\xi\right|^{1-\beta}} 1_{\{(\lambda_{1}+\xi)(\lambda_{2}-\xi)<0\}}d\xi
\right]W(d\lambda_{1})W(d\lambda_{2}),
\end{align}
and
\begin{align}\notag
Z^{(2)}(t)
=&8\int_{\mathbb{R}^{4}}^{'}e^{it(\lambda_{1}+\lambda_{2}+\lambda_{3}+\lambda_{4})}\frac{\widehat{\psi}(\lambda_{1}+\lambda_{2})
\widehat{\psi}(\lambda_{3}+\lambda_{4})}{|\lambda_{1}\lambda_{2}\lambda_{3}\lambda_{4}|^{(1-\beta)/2}}1_{\{\lambda_{1}\lambda_{2}<0\}}1_{\{\lambda_{3}\lambda_{4}<0\}}
\\\label{def:Zlimit2}&\times
1_{\{(\lambda_{1}+\lambda_{2})(\lambda_{3}+\lambda_{4})<0\}}
W(d\lambda_{1})W(d\lambda_{2})W(d\lambda_{3})W(d\lambda_{4}).
\end{align}
\end{Theorem}

\begin{Remark}
Because of Assumption \ref{Assumption:spectral}, which states that $g_{\tau}(s)\rightarrow 0$ as $s\rightarrow\infty$, (\ref{decay_in_g(s)}) together with the Cramer-Wold device implies that $\{Z_{s}(t)\}_{t\in \mathbb{R}}$ converges to $\{Z(t)\}_{t\in \mathbb{R}}$
in the finite-dimensional distribution sense as $s\rightarrow\infty$.
On the other hand, in addition to the condition $-\tau<(1-2\beta)/4$, if $\tau$ is strictly negative,
(\ref{decay_in_g(s)}) implies that
\begin{align}\notag
d_{\textup{Kol}}\left(\overset{n}{\underset{k=1}{\sum}}c_{k}Z_{s}\left(t_{k}\right),\overset{n}{\underset{k=1}{\sum}}c_{k}Z\left(t_{k}\right)\right)
\lesssim s^{2\tau/9+\delta}
\end{align}
for any $\delta>0$ because Theorem 1.5.3 in \cite{bingham1989regular} shows that
$
g_{\tau}(s)\lesssim s^{\tau+\delta}
$
for any $\delta>0$.
\end{Remark}
The proof of Theorem \ref{thm:nonStein} is detailed in Section \ref{sec:proof:thm:nonStein}. It mainly relies on the analytical probability method \cite{MR4003569}. Additionally, the proof involves the utilization of the following lemma concerning singular integrals and product contractions,
whose proof is provided in Section \ref{sec:proof:lemma:integral:nu12}.
\begin{Lemma}\label{lemma:integral:nu12}
For $\nu_{1},\nu_{2},\nu_{3}\in(-\beta,1-\beta)$, $\nu_{1}+\nu_{3}<1-2\beta$,
and $\nu_{2}+\nu_{3}<1-2\beta$,
\begin{align*}
\int_{\mathbb{R}^{2}}
\left|\lambda_{1}\right|^{\nu_{1}+\beta-1}\left|\lambda_{2}\right|^{\nu_{2}+\beta-1}
\left|\int_{\mathbb{R}}
\frac{|\widehat{\psi}(\lambda_{1}+\xi)\widehat{\psi}(\lambda_{2}-\xi)|}{\left|\xi\right|^{1-\beta-\nu_{3}}} d\xi
\right|^{2}
d\lambda_{1}d\lambda_{2}<\infty.
\end{align*}
\end{Lemma}

\begin{Remark}
About the distance between the $n$-dimensional random vector $\mathbf{Z}_{s}:=(Z_{s}(t_{1}),Z_{s}(t_{2}),...,Z_{s}(t_{n}))$ and
its limit $\mathbf{Z}:=(Z(t_{1}),Z(t_{2}),...,Z(t_{n}))$, the more clear approach is to consider their multidimensional Kolmogorov distance, denoted as $K(\mathbf{Z}_{s},\mathbf{Z})$.
An interesting open problem is to obtain the rate of convergence to zero of $K(\mathbf{Z}_{s},\mathbf{Z})$.
For the multivariate central limit theorems derived from the normalized least squares estimators of the regression coefficients with long-range dependence errors
and the rescaled solution of Burger's equation with random initial data,
the rates of convergence of the multidimensional Kolmogorov distance were obtained
by introduced normalizing
matrix-valued functions \cite{MR1642664,leonenko2003rate}.
Because the limiting process $Z$ in Theorem 1 is comprised of the second and fourth Wiener chaos,
the application of the methods in \cite{MR1642664,leonenko2003rate} to the multivariate non-central limit is not direct.
\end{Remark}

\subsection{Central limit theorems}

In the case where $2\beta>1$, we investigate the process
$\{s^{1/2}U[s]Y(st)\mid t\in \mathbb{R}\}$.
This process represents the normalized magnitude of the analytic wavelet transform of $Y$.
It's worth noting that the normalization considered in this case differs from that in (\ref{def:Zs}).
The following theorem will show that
its distribution will be close to the chi distribution with two degrees of freedom when the scale variable $s$ tends to infinity.
For technical reasons, we use the Wasserstein metric to measure the distance between them.

\begin{Definition}\label{def:W:distance}
Consider two random variables $V_{1}$ and $V_{2}$.
The Wasserstein distance between the distributions of $V_{1}$ and $V_{2}$, denoted as $d_{\textup{W}}(V_{1},V_{2})$, is defined as follows
\begin{align*}
d_{\textup{W}}(V_{1},V_{2}) = \underset{\|h\|_{\textup{Lip}}\leq 1}{\sup}
\left|\mathbb{E}[h(V_{1})]-\mathbb{E}[h(V_{2})]\right|,
\end{align*}
where the supremum is taken over functions with Lipschitz constants less than one.
\end{Definition}

\begin{Theorem}\label{thm:Stein}
Consider the random process $Y$ defined in (\ref{def:envelopeX}).
If the nonvanishing wavelet $\psi$ and its Fourier transform $\widehat{\psi}$ belong to $L^{1}$, and the assumption (\ref{condition_f}) about the spectral density function $f_{X}$ holds,
then for any $n\in \mathbb{N}$, distinct time points $t_{1},...,t_{n}\in \mathbb{R}$, and
$c_{1},...,c_{n}\in \mathbb{R}$
with $c_{1}^{2}+\cdots +c_{n}^{2}\leq 1$, when $s\rightarrow\infty$,
\begin{align}\notag
&d_{\textup{W}}\left(\overset{n}{\underset{k=1}{\sum}}c_{k}s^{1/2}U[s]Y(st_{k}),
\overset{n}{\underset{k=1}{\sum}}c_{k}\sqrt{N_{2k-1,s}^{2}+N_{2k,s}^{2}}\right)
\\\label{thm:nonStein:inequality}\lesssim &\left\{
\begin{array}{cl}
s^{-1/2} &\ \textup{for}\ \beta\in(3/4,1), \\
s^{-1/2+2(3/4-\beta)}L^{2}(s) &\ \textup{for}\ \beta\in(1/2,3/4).
\end{array}
\right.
\end{align}
Here, $(N_{1,s},N_{2,s},...,N_{2n-1,s},N_{2n,s})$ is a centered normal random vector and its covariance
matrix converges to $16\|f_{X}\|_{2}^{2}\Psi$ when $s\rightarrow\infty$.
The matrix $\mathrm{\Psi}=[\Psi(\ell,\ell')]_{1\leq\ell,\ell'\leq 2n}$ is positive definite
for any distinct time points $t_{1},...,t_{n}\in \mathbb{R}$ and composed of
\begin{equation}\label{def:Phi:bigmatrix}
\Psi(\ell,\ell') =
\left\{
\hspace{-0.2cm}\begin{array}{l}
\int_{0}^{\infty}\cos((t_{k}-t_{k'})u)|\widehat{\psi}(u)|^{2}du\  \textup{for}\ (\ell,\ell')=(2k-1,2k'-1)\ \textup{or}\ (2k,2k'),
\\
\int_{0}^{\infty}\sin((t_{k}-t_{k'})u)|\widehat{\psi}(u)|^{2}du\ \textup{for}\ (\ell,\ell')=(2k-1,2k'),
\\
\int_{0}^{\infty}\sin((t_{k'}-t_{k})u)|\widehat{\psi}(u)|^{2}du\ \textup{for}\ (\ell,\ell')=(2k,2k'-1),
\end{array}
\right.
\end{equation}
where $k,k'=1,2,...,n$.
\end{Theorem}

\begin{Remark}
By the continuous mapping theorem,
Theorem \ref{thm:Stein} shows that the process $\{sU^{2}[s]Y(st)\mid t\in \mathbb{R}\}$,
derived from the normalized scalogram of the process $Y$,
converges to a $\chi^{2}_{2}$-distributed random process in the finite-dimensional distribution sense
as $s\rightarrow\infty$.
It's noteworthy that the upper bounds of the Wasserstein distance, as shown in Theorem \ref{thm:Stein}, don't entirely align with the standard Berry-Esseen bound in the Wasserstein metric for the central limit theorem \cite{bonis2020stein}.
Specifically, Theorem \ref{thm:Stein} shows that the speed of convergence to zero decreases from $s^{-1/2}L^{2}(s)$ to $s^{0}L^{2}(s)$ according to the function $s^{-1/2+2(3/4-\beta)}L^{2}(s)$
when the long-memory parameter $\beta$ of the initial input $X$ decreases from 3/4 to 1/2.
This observation resonates with findings in the literature \cite{bierme2011central,breton2008error} regarding
the total variation distance between the normalized sum of the Hermite power of increments of a fractional Brownian motion and its limit,
as well as how the rate of convergence depends on the Hurst index of the fractional Brownian motion.
\end{Remark}

The proof of Theorem \ref{thm:Stein} relies on the multivariate Stein method and the Malliavin calculus \cite{nourdin2009stein,nourdin2012normal}.
There are two key points in this approach.
The first concerns the covariance matrices of
the normal random vectors used for the normal approximation.
Remark 5.1.4 and Theorem 6.1.1 in \cite{nourdin2009stein} highlight the challenge in ensuring the positive definiteness of the covariance matrices.
By using the observation (\ref{real_complex_link}) along with the characterization of strictly positive definite functions
in \cite[Corollary 6.9]{wendland2004scattered}, we have proven the positive definiteness of the covariance matrix for
$(N_{1,s},N_{2,s},...,N_{2n-1,s},N_{2n,s})$.
The second key point is about the estimate of the norms of product contractions (\ref{estimate:singularintegral1}) of the integrand functions of multiple Wiener integrals.
Lemma \ref{lemma:singular_integral} below plays a crucial role in estimating these norms.
Because Lemma \ref{lemma:singular_integral} has not appeared in existing literature, its proof is provided in Section \ref{sec:proof:lemma:singular_integral}.
\begin{Lemma}\label{lemma:singular_integral}
If the wavelet $\psi$ and its Fourier transform $\widehat{\psi}$ belong to $L^{1}$ and  $\beta\in(1/2,3/4)$,
then for any
$\mu_{1},\mu_{2},\mu_{3},\mu_{4}\in\{-1,1\}$, the integral below
\begin{align}\label{def:singular_integral}
\mathcal{K}:=\int_{\mathbb{R}^{4}}
 \frac{|\widehat{\psi}(x)| |\widehat{\psi}(y)||\widehat{\psi}(z)| |\widehat{\psi}(x+y-z)|}{
|u|^{1-\beta-\delta \mu_{1}}|x+y-u|^{1-\beta-\delta\mu_{2}}|x-u|^{1-\beta-\delta\mu_{3}}|z-u|^{1-\beta-\delta\mu_{4}}}
\,du\,dx\,dy\, dz
\end{align}
is finite for sufficiently small $\delta>0$.
\end{Lemma}

\begin{Corollary}\label{corollary:cov:chi}
Under the hypotheses of Theorem \ref{thm:Stein}, the covariance between any time points $t$ and $t'$
of the limiting process of
$sU^{2}[s]Y(s\cdot)$
can be expressed as follows
\begin{align*}
2^{10}\|f_{X}\|_{2}^{4}\left[ \left(\int_{0}^{\infty}\cos((t-t')u)|\widehat{\psi}(u)|^{2}du\right)^{2}
+\left(\int_{0}^{\infty}\sin((t-t')u)|\widehat{\psi}(u)|^{2}du\right)^{2}\right].
\end{align*}
\end{Corollary}

\section{Concluding remarks}\label{sec:Conclusion}
For the nonlinear transformation depicted in Figure \ref{FlowChart}, from Theorems \ref{thm:nonStein} and \ref{thm:Stein},  we observed that the long-memory parameter of the initial input $X$ at $1/2$ acts as a dividing point.
Remarkably, this division directly corresponds to the long-memory parameter of $Y$, which represents the squared modulus of the analytic counterpart of $X$.
When the long-memory parameter of $X$ either falls below or exceeds this threshold, it determines the nature of the process $Y$ -- whether it exhibits long-range dependence or short-range dependence.
This distinction leads to the convergence of the rescaled version of $U[s]Y$ towards different limits.
The former (i.e., $s^{2\beta}L^{-2}(s)U^{2}[s]Y(s\cdot)$) converges to an additive combination of the second and fourth Wiener chaos, while the latter (i.e., $sU^{2}[s]Y(s\cdot)$) tends towards a chi-square random process.
Besides, we also conducted a systematic analysis of
the relationship between the convergence rate of the large scaling limits and the strength of long-range dependence in the initial input $X$.
About the upper bound of the distance between the distributions of the rescaled $U[s]Y$ and its limit, our observations indicate that the speed of convergence to zero diminishes as the strength of long-range dependence in the initial input $X$ increases.


\section{Proofs}\label{sec:proof}

\subsection{Proof of Lemma \ref{lemma:chisquare}}\label{sec:proof:lemma:chisquare}
First of all, $\{Y(t)\}_{t\in \mathbb{R}}$ can be rewritten as  $Y(t)= |X(t)+i\widetilde{X}(t)|^{2}$,
where
\begin{align*}
\widetilde{X}(t) = -i\int_{\mathbb{R}}
e^{i\lambda t}\textup{sgn}(\lambda)\sqrt{f_{X}(\lambda)}W(d\lambda).
\end{align*}
From (\ref{ortho}), we know that
$\widetilde{X}$ is also a real-valued
Gaussian process and
\begin{align}\label{cov:NRNI}
\mathbb{E}\left[\begin{array}{cc}X(t)X(t')& X(t)\widetilde{X}(t')\\ \widetilde{X}(t)X(t') & \widetilde{X}(t)\widetilde{X}(t')\end{array}\right] =
2\int_{0}^{\infty}\left[\begin{array}{lr}\cos(\lambda(t-t^{'})) & -\sin(\lambda(t-t^{'}))\\ \sin(\lambda(t-t'))  &\cos(\lambda(t-t^{'})) \end{array}\right]f_{X}(\lambda)d\lambda.
\end{align}
Especially, $\mathbb{E}[X^{2}(t)] = \mathbb{E}[\widetilde{X}^{2}(t)] = \|f_{X}\|_{1}$ and
$\mathbb{E}[X(t)\widetilde{X}(t)]=0$, which implies the independence between $X(t)$ and $\widetilde{X}(t)$.
Hence, $Y(t)/\|f_{X}\|_{1}=(X^{2}(t)+\widetilde{X}^{2}(t))/\|f_{X}\|_{1}$ follows a chi-square distribution with two degrees of freedom.

By the product formula \cite{major1981lecture},
\begin{align}\label{xsquare}
X^{2}(t) = \|f_{X}\|_{1}+\int_{\mathbb{R}^{2}}^{'}
e^{i(\lambda_{1}+\lambda_{2}) t}\sqrt{f_{X}(\lambda_{1})f_{X}(\lambda_{2})}W(d\lambda_{1})W(d\lambda_{2})
\end{align}
and
\begin{align}\label{xtildesquare}
\widetilde{X}^{2}(t) = \|f_{X}\|_{1}-\int_{\mathbb{R}^{2}}^{'}
e^{i(\lambda_{1}+\lambda_{2}) t}\textup{sgn}(\lambda_{1}\lambda_{2})\sqrt{f_{X}(\lambda_{1})f_{X}(\lambda_{2})}W(d\lambda_{1})W(d\lambda_{2}).
\end{align}
The representation (\ref{Ychaosdep})  is obtained by combining (\ref{xsquare}) and (\ref{xtildesquare}).
The covariance function of $Y$
follows from (\ref{cov:NRNI}). More precisely,
\begin{align}\notag
&\textup{Cov}(Y(t),Y(t')) = \textup{Cov}(X^{2}(t)+\widetilde{X}^{2}(t),X^{2}(t')+\widetilde{X}^{2}(t'))
\\\notag=&2\left\{[\textup{Cov}(X(t),X(t'))]^{2}
+[\textup{Cov}(X(t),\widetilde{X}(t'))]^{2}
+[\textup{Cov}(\widetilde{X}(t),X(t'))]^{2}
+[\textup{Cov}(\widetilde{X}(t),\widetilde{X}(t'))]^{2}\right\}
\\\notag=&16\left|\int_{0}^{\infty} e^{i\lambda(t-t')}f_{X}(\lambda)d\lambda\right|^{2}.
\end{align}

\subsection{Proof of Lemma \ref{lemma:chasodep:UY}}\label{sec:proof:lemma:chasodep:UY}
By the spectral representation of $Y$ in (\ref{Ychaosdep}), the stochastic Fubini theorem \cite[Theorem 2.1]{pipiras2010regularization},
and $\widehat{\psi_{s}}(\lambda)=\widehat{\psi}(s\lambda)$,
\begin{align}\label{spectral_Xstarpsi}
Y\star \psi_{s}(t) = \int_{\mathbb{R}^{2}}^{'}e^{it(\lambda_{1}+\lambda_{2})} k(\lambda_{1},\lambda_{2})W(d\lambda_{1})W(d\lambda_{2}),
\end{align}
where the influence of $\mathbb{E}[Y]$ on $Y\star \psi_{s}$ has been eliminated by the property $\int_{\mathbb{R}}\psi(u)du = 0$,
and
\begin{align*}
k(\lambda_{1},\lambda_{2}) =2\sqrt{f_{X}(\lambda_{1})
f_{X}(\lambda_{2})}1_{\{\lambda_{1}\lambda_{2}<0\}}\widehat{\psi}(s(\lambda_{1}+\lambda_{2})).
\end{align*}
By the product formula \cite{major1981lecture},
\begin{align}\notag
\left(Y\star \psi_{s}(t)\right)^{2}
=& 2\|k\|_{2}^{2}+
4\int_{\mathbb{R}^{2}}^{'}
e^{it(\lambda_{1}+\lambda_{2})}k\otimes_{1}k(\lambda_{1},\lambda_{2})
W(d\lambda_{1})W(d\lambda_{2})
\\\label{spectral_U[s,s']comp1}&+\int_{\mathbb{R}^{4}}^{'}
e^{it(\lambda_{1}+\lambda_{2}+\lambda_{3}+\lambda_{4})}k\otimes k(\lambda_{1},\lambda_{2},\lambda_{3},\lambda_{4})
\overset{4}{\underset{\ell=1}{\prod}}W(d\lambda_{\ell}),
\end{align}
where
\begin{align*}
k\otimes_{1}k(\lambda_{1},\lambda_{2}) = \int_{\mathbb{R}}k(\lambda_{1},\xi)k(\lambda_{2},-\xi)d\xi
\end{align*}
and
$k\otimes k(\lambda_{1},\lambda_{2},\lambda_{3},\lambda_{4}) =k(\lambda_{1},\lambda_{2})k(\lambda_{3},\lambda_{4})$.

By the convolution theorem, $\widehat{\mathcal{H}\psi_{s}}(\lambda)=-i\textup{sgn}(\lambda)\widehat{\psi}(s\lambda)$, where
$\textup{sgn}(\lambda)=\lambda/|\lambda|$ for $\lambda\neq0$ and $\textup{sgn}(0)=0$.
In parallel with (\ref{spectral_Xstarpsi})
and (\ref{spectral_U[s,s']comp1}),
\begin{align*}
Y\star \mathcal{H}\psi_{s}(t)=-i \int_{\mathbb{R}^{2}}^{'}e^{it(\lambda_{1}+\lambda_{2})}
\textup{sgn}(\lambda_{1}+\lambda_{2})
k(\lambda_{1},\lambda_{2}) W(d\lambda_{1})W(d\lambda_{2})
\end{align*}
and
\begin{align}\notag
\left(Y\star \mathcal{H}\psi_{s}(t)\right)^{2}
=& 2\|k\|_{2}^{2}+
4\int_{\mathbb{R}^{2}}^{'}
e^{it(\lambda_{1}+\lambda_{2})}\widetilde{k}\otimes_{1}\widetilde{k}(\lambda_{1},\lambda_{2})
W(d\lambda_{1})W(d\lambda_{2})
\\\label{spectral_U[s,s']comp2}&+\int_{\mathbb{R}^{4}}^{'}
e^{it(\lambda_{1}+\lambda_{2}+\lambda_{3}+\lambda_{4})}\widetilde{k}\otimes \widetilde{k}(\lambda_{1},\lambda_{2},\lambda_{3},\lambda_{4})
\overset{4}{\underset{\ell=1}{\prod}}W(d\lambda_{\ell}),
\end{align}
where $\widetilde{k}(\lambda_{1},\lambda_{2})=
-i\textup{sgn}(\lambda_{1}+\lambda_{2})k(\lambda_{1},\lambda_{2}).$
The proof is finished by combining (\ref{spectral_U[s,s']comp1}) and (\ref{spectral_U[s,s']comp2}).

\subsection{Proof of Theorem \ref{thm:nonStein}}\label{sec:proof:thm:nonStein}
Here, we focus on the case where $n=1$, $c_{1}=1$, and $t_{1}=t$ to simplify notation. The following proof can be extended to the general case for $n\geq2.$
First of all, we denote $Z_{s}$ as the sum of three components: an $s$-dependent constant $Z_s^{(0)}$ and two random processes, $Z_s^{(1)}$ and $Z_s^{(2)}$, where
\begin{align}\label{def:Z0s}
Z_{s}^{(0)}=&s^{2\beta}L^{-2}(s)\mathbb{E}[U^{2}[s]Y(st)],
\end{align}
\begin{align}\notag
Z_{s}^{(1)}(t)=s^{2\beta}L^{-2}(s)\int_{\mathbb{R}^{2}}^{'}e^{ist(\lambda_{1}+\lambda_{2})}
K_{2}(\lambda_{1},\lambda_{2})
W(d\lambda_{1})W(d\lambda_{2}),
\end{align}
and
\begin{align}\notag
Z_{s}^{(2)}(t) = s^{2\beta}L^{-2}(s)\int_{\mathbb{R}^{4}}^{'}e^{ist(\lambda_{1}+\lambda_{2}+\lambda_{3}+\lambda_{4})}
K_{4}(\lambda_{1},\lambda_{2},\lambda_{3},\lambda_{4})
\overset{4}{\underset{\ell=1}{\prod}}W(d\lambda_{\ell}),
\end{align}
where $K_{2}$ and $K_{4}$ are defined in Lemma \ref{lemma:chasodep:UY}.
Consider two random processes  $V_{1}$ and $V_{2}$, and denote $V_{1}\overset{d}{=}V_{2}$ to signify that their finite-dimensional distributions are identical.
Utilizing the self-similarity property $W(s^{-1}d\lambda)\overset{d}{=}s^{-1/2}W(d\lambda)$ (see also \cite[Lemma 6.1]{taqqu1979convergence}), we have
\begin{align}\label{dist_equiv}
Z_{s}(t)=Z^{(0)}_{s}+Z^{(1)}_{s}(t)
+ Z^{(2)}_{s}(t)
\overset{d}{=}\widetilde{Z}_{s}(t),
\end{align}
where $\widetilde{Z}_{s}(t)=Z^{(0)}_{s}+\widetilde{Z}^{(1)}_{s}(t)+\widetilde{Z}^{(2)}_{s}(t)$ with
\begin{align}\notag
\widetilde{Z}^{(1)}_{s}(t) =&
s^{2\beta-1}L^{-2}(s)\int_{\mathbb{R}^{2}}^{'}e^{it(\lambda_{1}+\lambda_{2})}
K_{2}(\lambda_{1}/s,\lambda_{2}/s)
W(d\lambda_{1})W(d\lambda_{2})
\\\notag=&
8s^{2\beta-2}L^{-2}(s)\int_{\mathbb{R}^{2}}^{'}e^{it(\lambda_{1}+\lambda_{2})}
\sqrt{f_{X}\left(\frac{\lambda_{1}}{s}\right)
f_{X}\left(\frac{\lambda_{2}}{s}\right)}1_{\{\lambda_{1}\lambda_{2}<0\}}
\\\notag&\times\left[\int_{\mathbb{R}}f_{X}\left(\frac{\xi}{s}\right)
\widehat{\psi}(\lambda_{1}+\xi)\widehat{\psi}(\lambda_{2}-\xi) 1_{\{(\lambda_{1}+\xi)(\lambda_{2}-\xi)<0\}}d\xi\right]
W(d\lambda_{1})W(d\lambda_{2})
\\\notag=&
8s^{2\beta-2}L^{-2}(s)\int_{\mathbb{R}^{2}}^{'}e^{it(\lambda_{1}+\lambda_{2})}
L^{1/2}\left(\frac{s}{|\lambda_{1}|}\right)L^{1/2}\left(\frac{s}{|\lambda_{2}|}\right)
\left|\frac{\lambda_{1}}{s}\right|^{\frac{\beta-1}{2}}
\left|\frac{\lambda_{2}}{s}\right|^{\frac{\beta-1}{2}}
1_{\{\lambda_{1}\lambda_{2}<0\}}
\\\notag&\times\left[\int_{\mathbb{R}}
L\left(\frac{s}{|\xi|}\right)
\left|\frac{\xi}{s}\right|^{\beta-1}
\widehat{\psi}(\lambda_{1}+\xi)\widehat{\psi}(\lambda_{2}-\xi) 1_{\{(\lambda_{1}+\xi)(\lambda_{2}-\xi)<0\}}d\xi\right]
W(d\lambda_{1})W(d\lambda_{2})
\\\label{def:Z1equiv}=&
8\int_{\mathbb{R}^{2}}^{'}e^{it(\lambda_{1}+\lambda_{2})}
\frac{L^{1/2}\left(s/|\lambda_{1}|\right)}{L^{1/2}(s)}
\frac{L^{1/2}\left(s/|\lambda_{2}|\right)}{L^{1/2}(s)}
\left|\lambda_{1}\lambda_{2}\right|^{\frac{\beta-1}{2}}
1_{\{\lambda_{1}\lambda_{2}<0\}}
\\\notag&\times\left[\int_{\mathbb{R}}
\frac{L\left(s/|\xi|\right)}{L(s)}
\frac{\widehat{\psi}(\lambda_{1}+\xi)\widehat{\psi}(\lambda_{2}-\xi)}{\left|\xi\right|^{1-\beta}} 1_{\{(\lambda_{1}+\xi)(\lambda_{2}-\xi)<0\}}d\xi\right]
W(d\lambda_{1})W(d\lambda_{2})
\end{align}
and
\begin{align}\notag
\widetilde{Z}_{s}^{(2)}(t) =& s^{2(\beta-1)}L^{-2}(s)\int_{\mathbb{R}^{4}}^{'}e^{it(\lambda_{1}+\cdots+\lambda_{4})}
K_{4}(\lambda_{1}/s,...,\lambda_{4}/s)
W(d\lambda_{1})W(d\lambda_{2})W(d\lambda_{3})W(d\lambda_{4})
\\\notag=&8\int_{\mathbb{R}^{4}}^{'}e^{it(\lambda_{1}+\cdots+\lambda_{4})}
\left[\overset{4}{\underset{\ell=1}{\prod}}\frac{L^{1/2}(s/\lambda_{\ell})}{L^{1/2}(s)}\right]
1_{\{\lambda_{1}\lambda_{2}<0\}}1_{\{\lambda_{3}\lambda_{4}<0\}}
1_{\{(\lambda_{1}+\lambda_{2})(\lambda_{3}+\lambda_{4})<0\}}
\\\label{def:Z2equiv}&\times\frac{\widehat{\psi}(\lambda_{1}+\lambda_{2})
\widehat{\psi}(\lambda_{3}+\lambda_{4})}{|\lambda_{1}\lambda_{2}\lambda_{3}\lambda_{4}|^{(1-\beta)/2}}
W(d\lambda_{1})W(d\lambda_{2})W(d\lambda_{3})W(d\lambda_{4}).
\end{align}
By (\ref{dist_equiv}) and \cite[Lemma 1.8]{MR1353441}, for any $\varepsilon>0$,
\begin{align}\notag
&d_{\textup{Kol}}\left(Z_{s}(t),Z(t)\right)
= d_{\textup{Kol}}\left(\widetilde{Z}_{s}(t),Z(t)\right)
\\\notag=&d_{\textup{Kol}}\left(Z(t)+(Z^{(0)}_{s}-Z^{(0)})+(\widetilde{Z}_{s}^{(1)}(t)-Z^{(1)}(t))+(\widetilde{Z}_{s}^{(2)}(t)-Z^{(2)}(t)),Z(t)\right)
\\\notag\leq& \overset{2}{\underset{\ell=1}{\sum}} \mathbb{P}\left(|\widetilde{Z}^{(\ell)}_{s}(t)-Z^{(\ell)}(t)|\geq \frac{\varepsilon}{2}\right)+d_{\textup{Kol}}\left(Z(t)+Z^{(0)}_{s}-Z^{(0)},Z(t)\right)
\\\notag&+d_{\textup{Kol}}\left(Z(t)+\varepsilon,Z(t)\right)
\\\notag\leq& \frac{4}{\varepsilon^{2}}\overset{2}{\underset{\ell=1}{\sum}}\mathbb{E}[|\widetilde{Z}^{(\ell)}_{s}(t)-Z^{(\ell)}(t)|^{2}]
+d_{\textup{Kol}}\left(Z(t)+Z^{(0)}_{s}-Z^{(0)},Z(t)\right)
\\\label{kol:distance1}&+d_{\textup{Kol}}\left(Z(t)+\varepsilon,Z(t)\right).
\end{align}

\noindent{$\bullet$ \it Upper bound of $\mathbb{E}[|\widetilde{Z}^{(1)}_{s}(t)-Z^{(1)}(t)|^{2}]$}:
By the definitions of $Z^{(1)}_{s}$ and $\widetilde{Z}^{(1)}$ in (\ref{def:Zlimit1}) and (\ref{def:Z1equiv}),
\begin{align}\notag
\widetilde{Z}^{(1)}_{s}(t)-Z^{(1)}(t)
=&
8\int_{\mathbb{R}^{2}}^{'}e^{it(\lambda_{1}+\lambda_{2})}
\left|\lambda_{1}\lambda_{2}\right|^{\frac{\beta-1}{2}}
\\\notag&\times
\int_{\mathbb{R}}\hspace{-0.1cm}
D_{1}(\lambda_{1},\lambda_{2},\xi)\frac{\widehat{\psi}(\lambda_{1}+\xi)\widehat{\psi}(\lambda_{2}-\xi)}{\left|\xi\right|^{1-\beta}} d\xi\,
W(d\lambda_{1})W(d\lambda_{2}),
\end{align}
where
\begin{align}\notag
D_{1}(\lambda_{1},\lambda_{2},\xi) =
\left[\left(\frac{L\left(s/|\lambda_{1}|\right)}{L(s)}\right)^{1/2}\left(\frac{L\left(s/|\lambda_{2}|\right)}{L(s)}\right)^{1/2}
\frac{L\left(s/|\xi|\right)}{L(s)}-1\right]
1_{\{\lambda_{1}\lambda_{2}<0\}}1_{\{(\lambda_{1}+\xi)(\lambda_{2}-\xi)<0\}}.
\end{align}
By the isometry property of multiple Wiener integrals,
\begin{align}\notag
&\mathbb{E}\left[|\widetilde{Z}^{(1)}_{s}(t)-Z^{(1)}(t)|^{2}\right]
\\\label{ineq:L2:Zs-Z:v1}=&
128\int_{\mathbb{R}^{2}}
\left|\lambda_{1}\lambda_{2}\right|^{\beta-1}
\left|\int_{\mathbb{R}}
D_{1}(\lambda_{1},\lambda_{2},\xi)\frac{\widehat{\psi}(\lambda_{1}+\xi)\widehat{\psi}(\lambda_{2}-\xi)}{\left|\xi\right|^{1-\beta}} d\xi
\right|^{2}
d\lambda_{1}d\lambda_{2}.
\end{align}
Applying the inequality $|xy-1|\leq |x^{2}-1|+|y^{2}-1|$ for $x, y \geq 0$ repeatedly yields
\begin{align}\notag
|D_{1}(\lambda_{1},\lambda_{2},\xi)|\leq&
\left|\left(\frac{L\left(s/|\lambda_{1}|\right)}{L(s)}\right)^{1/2}\left(\frac{L\left(s/|\lambda_{2}|\right)}{L(s)}\right)^{1/2}
\frac{L\left(s/|\xi|\right)}{L(s)}-1\right|
\\\label{ineq:D}\leq& \left|\frac{L^{2}\left(s/|\lambda_{1}|\right)}{L^{2}(s)}-1\right|+\left|\frac{L^{2}\left(s/|\lambda_{2}|\right)}{L^{2}(s)}-1\right|
+\left|\frac{L^{2}\left(s/|\xi|\right)}{L^{2}(s)}-1\right|.
\end{align}
By (\ref{ineq:L2:Zs-Z:v1}) and  (\ref{ineq:D}),
\begin{align*}
\mathbb{E}\left[|\widetilde{Z}^{(1)}_{s}(t)-Z^{(1)}(t)|^{2}\right]
\lesssim
I_{1}+I_{2},
\end{align*}
where
\begin{align*}
I_{1} =
\int_{\mathbb{R}^{2}}
\left|\lambda_{1}\lambda_{2}\right|^{\beta-1}\left|\frac{L^{2}\left(s/|\lambda_{1}|\right)}{L^{2}(s)}-1\right|^{2}
\left[\int_{\mathbb{R}}
\frac{|\widehat{\psi}(\lambda_{1}+\xi)\widehat{\psi}(\lambda_{2}-\xi)|}{\left|\xi\right|^{1-\beta}} d\xi
\right]^{2}
d\lambda_{1}d\lambda_{2}
\end{align*}
and
\begin{align*}
I_{2} =
\int_{\mathbb{R}^{2}}
\left|\lambda_{1}\lambda_{2}\right|^{\beta-1}
\left[\int_{\mathbb{R}}
\left|\frac{L^{2}\left(s/|\xi|\right)}{L^{2}(s)}-1\right|
\frac{|\widehat{\psi}(\lambda_{1}+\xi)\widehat{\psi}(\lambda_{2}-\xi)|}{\left|\xi\right|^{1-\beta}} d\xi
\right]^{2}
d\lambda_{1}d\lambda_{2}.
\end{align*}

Let $B$ be the interval $[-1,1]$, and let $B^c$ denote its complement.
Then, $I_{1} = I_{1,1}+I_{1,2}$, where
\begin{align*}
I_{1,1} =
\int_{\mathbb{R}}\int_{B}
\left|\lambda_{1}\lambda_{2}\right|^{\beta-1}\left|\frac{L^{2}\left(s/|\lambda_{1}|\right)}{L^{2}(s)}-1\right|^{2}
\left[\int_{\mathbb{R}}
\frac{|\widehat{\psi}(\lambda_{1}+\xi)\widehat{\psi}(\lambda_{2}-\xi)|}{\left|\xi\right|^{1-\beta}} d\xi
\right]^{2}
d\lambda_{1}d\lambda_{2}
\end{align*}
and
\begin{align*}
I_{1,2} =&
\int_{\mathbb{R}}\int_{B^{c}}
\left|\lambda_{1}\lambda_{2}\right|^{\beta-1}\left|\frac{L^{2}\left(s/|\lambda_{1}|\right)}{L^{2}(s)}-1\right|^{2}
\left[\int_{\mathbb{R}}
\frac{|\widehat{\psi}(\lambda_{1}+\xi)\widehat{\psi}(\lambda_{2}-\xi)|}{\left|\xi\right|^{1-\beta}} d\xi
\right]^{2}
d\lambda_{1}d\lambda_{2}.
\end{align*}
By Lemma \ref{lemma:L:ratio:easy} in  \ref{sec:appendix:L},
\begin{align*}
I_{1,1}\lesssim
g_{\tau}^{2}(s)\int_{\mathbb{R}}\int_{B}
\left|\lambda_{1}\right|^{\beta-1-2\delta}\left|\lambda_{2}\right|^{\beta-1}
\left[\int_{\mathbb{R}}
\frac{|\widehat{\psi}(\lambda_{1}+\xi)\widehat{\psi}(\lambda_{2}-\xi)|}{\left|\xi\right|^{1-\beta}} d\xi
\right]^{2}
d\lambda_{1}d\lambda_{2}
\end{align*}
for any $\delta>0$. The integral above is finite if $2\beta<1$ and $\delta$ is sufficiently small, as demonstrated by Lemma \ref{lemma:integral:nu12}.
By Lemma \ref{lemma:L:ratio:diff} in  \ref{sec:appendix:L},
\begin{align*}
I_{1,2}\lesssim  g_{\tau}^{2}(s)\int_{\mathbb{R}}\int_{B^{c}}
\left|\lambda_{1}\right|^{\beta-1-4\tau+2\delta}\left|\lambda_{2}\right|^{\beta-1}
\left[\int_{\mathbb{R}}
\frac{|\widehat{\psi}(\lambda_{1}+\xi)\widehat{\psi}(\lambda_{2}-\xi)|}{\left|\xi\right|^{1-\beta}} d\xi
\right]^{2}
d\lambda_{1}d\lambda_{2}
\end{align*}
for any $\delta>0$. The integral above is finite if $-\tau<(1-2\beta)/4$ and $\delta$ is sufficiently small, as demonstrated by Lemma \ref{lemma:integral:nu12}.
Therefore, under the condition $-\tau<(1-2\beta)/4$,
\begin{align}\label{final:I1}
I_{1}\lesssim  g_{\tau}^{2}(s).
\end{align}

Observe that $I_{2} \leq 2(I_{2,1} + I_{2,2}),$ where
\begin{align*}
I_{2,1} =
\int_{\mathbb{R}^{2}}
\left|\lambda_{1}\lambda_{2}\right|^{\beta-1}
\left[\int_{B}
\left|\frac{L^{2}\left(s/|\xi|\right)}{L^{2}(s)}-1\right|
\frac{|\widehat{\psi}(\lambda_{1}+\xi)\widehat{\psi}(\lambda_{2}-\xi)|}{\left|\xi\right|^{1-\beta}} d\xi
\right]^{2}
d\lambda_{1}d\lambda_{2}
\end{align*}
and
\begin{align*}
I_{2,2} =
\int_{\mathbb{R}^{2}}
\left|\lambda_{1}\lambda_{2}\right|^{\beta-1}
\left[\int_{B^{c}}
\left|\frac{L^{2}\left(s/|\xi|\right)}{L^{2}(s)}-1\right|
\frac{|\widehat{\psi}(\lambda_{1}+\xi)\widehat{\psi}(\lambda_{2}-\xi)|}{\left|\xi\right|^{1-\beta}} d\xi
\right]^{2}
d\lambda_{1}d\lambda_{2}.
\end{align*}
By Lemma \ref{lemma:L:ratio:easy},
\begin{align}\label{estimate:J1:v2}
I_{2,1} \lesssim
g_{\tau}^{2}(s)\int_{\mathbb{R}^{2}}
\left|\lambda_{1}\lambda_{2}\right|^{\beta-1}
\left[\int_{B}
\frac{|\widehat{\psi}(\lambda_{1}+\xi)\widehat{\psi}(\lambda_{2}-\xi)|}{\left|\xi\right|^{1-\beta+\delta}} d\xi
\right]^{2}
d\lambda_{1}d\lambda_{2}
\end{align}
for any $\delta>0$.
Referring to Lemma \ref{lemma:integral:nu12}, the integral in (\ref{estimate:J1:v2}) is finite when $2\beta<1$ and $\delta$ is sufficiently small.
On the other hand, by Lemma \ref{lemma:L:ratio:diff},
\begin{align*}
I_{2,2}\lesssim g_{\tau}^{2}(s)
\int_{\mathbb{R}^{2}}
\left|\lambda_{1}\lambda_{2}\right|^{\beta-1}
\left[\int_{B^{c}}
\frac{|\widehat{\psi}(\lambda_{1}+\xi)\widehat{\psi}(\lambda_{2}-\xi)|}{\left|\xi\right|^{1-\beta+2\tau-\delta}} d\xi
\right]^{2}
d\lambda_{1}d\lambda_{2}
\end{align*}
for any $\delta>0$, where the integral is finite when $-\tau<(1-2\beta)/2$ and $\delta$ is sufficiently small, as shown in Lemma \ref{lemma:integral:nu12}.
Therefore,
\begin{equation}\label{final:I2}
I_{2} \lesssim g_{\tau}^{2}(s)
\end{equation}
under the condition $-\tau<(1-2\beta)/2$.

Combining (\ref{final:I1}) and (\ref{final:I2}) yields
\begin{align}\label{kol:distance2}
\mathbb{E}\left[|\widetilde{Z}^{(1)}(t)-Z^{(1)}(t)|^{2}\right]\lesssim  g_{\tau}^{2}(s)
\end{align}
under the condition $-\tau<(1-2\beta)/4$.
\\

\noindent{$\bullet$ \it Upper bound of $\mathbb{E}[|\widetilde{Z}^{(2)}_{s}(t)-Z^{(2)}(t)|^{2}]$}:
By the definitions of $Z^{(2)}$ and $\widetilde{Z}^{(2)}_{s}$ in (\ref{def:Zlimit2}) and  (\ref{def:Z2equiv}),
\begin{align}\notag
\widetilde{Z}^{(2)}_{s}(t)-Z^{(2)}(t)
=&8\int_{\mathbb{R}^{4}}^{'}e^{it(\lambda_{1}+\cdots+\lambda_{4})}
D_{2}(\lambda_{1},...,\lambda_{4})
\\\notag&\times
\frac{\widehat{\psi}(\lambda_{1}+\lambda_{2})
\widehat{\psi}(\lambda_{3}+\lambda_{4})}{|\lambda_{1}\lambda_{2}\lambda_{3}\lambda_{4}|^{(1-\beta)/2}}
W(d\lambda_{1})W(d\lambda_{2})W(d\lambda_{3})W(d\lambda_{4}),
\end{align}
where
\begin{align*}
D_{2}(\lambda_{1},\lambda_{2},\lambda_{3},\lambda_{4})
=\left[\overset{4}{\underset{\ell=1}{\prod}}\frac{L^{1/2}(s/\lambda_{\ell})}{L^{1/2}(s)}-1\right]
1_{\{\lambda_{1}\lambda_{2}<0\}}1_{\{\lambda_{3}\lambda_{4}<0\}}
1_{\{(\lambda_{1}+\lambda_{2})(\lambda_{3}+\lambda_{4})<0\}}.
\end{align*}
By the isometry property of multiple Wiener integrals and the Minkowski inequality,
\begin{align}\notag
\mathbb{E}\left[|\widetilde{Z}^{(2)}_{s}(t)-Z^{(2)}(t)|^{2}\right]
\leq&64\cdot4!\int_{\mathbb{R}^{4}}
|D_{2}(\lambda_{1},...,\lambda_{4})|^{2}
\\\label{D2appear}&\times\frac{|\widehat{\psi}(\lambda_{1}+\lambda_{2})
\widehat{\psi}(\lambda_{3}+\lambda_{4})|^{2}}{|\lambda_{1}\lambda_{2}\lambda_{3}\lambda_{4}|^{1-\beta}}
d\lambda_{1}d\lambda_{2}d\lambda_{3}d\lambda_{4}.
\end{align}
Applying the inequality
\begin{align*}
|D_{2}(\lambda_{1},...,\lambda_{4})|
\leq \overset{4}{\underset{\ell=1}{\sum}}
\left|\frac{L^{2}\left(s/|\lambda_{\ell}|\right)}{L^{2}(s)}-1\right|
\end{align*}
to (\ref{D2appear})
yields
\begin{align}\notag
\mathbb{E}\left[|\widetilde{Z}^{(2)}_{s}(t)-Z^{(2)}(t)|^{2}\right]
\lesssim \int_{\mathbb{R}^{4}}&
\left|\frac{L^{2}\left(s/|\lambda_{1}|\right)}{L^{2}(s)}-1\right|^{2}
\\\label{sumJ1J2}&\times\frac{|\widehat{\psi}(\lambda_{1}+\lambda_{2})
\widehat{\psi}(\lambda_{3}+\lambda_{4})|^{2}}{|\lambda_{1}\lambda_{2}\lambda_{3}\lambda_{4}|^{1-\beta}}
d\lambda_{1}d\lambda_{2}d\lambda_{3}d\lambda_{4}.
\end{align}
Denote
\begin{align}\notag
J_{1}= \int_{\mathbb{R}^{3}}\int_{B}
\left|\frac{L^{2}\left(s/|\lambda_{1}|\right)}{L^{2}(s)}-1\right|^{2}
\frac{|\widehat{\psi}(\lambda_{1}+\lambda_{2})
\widehat{\psi}(\lambda_{3}+\lambda_{4})|^{2}}{|\lambda_{1}\lambda_{2}\lambda_{3}\lambda_{4}|^{1-\beta}}
d\lambda_{1}d\lambda_{2}d\lambda_{3}d\lambda_{4},
\end{align}
where $B$ is the interval $[-1,1]$, and
\begin{align}\notag
J_{2}= \int_{\mathbb{R}^{3}}\int_{B^{c}}
\left|\frac{L^{2}\left(s/|\lambda_{1}|\right)}{L^{2}(s)}-1\right|^{2}
\frac{|\widehat{\psi}(\lambda_{1}+\lambda_{2})
\widehat{\psi}(\lambda_{3}+\lambda_{4})|^{2}}{|\lambda_{1}\lambda_{2}\lambda_{3}\lambda_{4}|^{1-\beta}}
d\lambda_{1}d\lambda_{2}d\lambda_{3}d\lambda_{4}.
\end{align}
By Lemma \ref{lemma:L:ratio:easy}, for sufficiently small $\delta>0$,
\begin{align}\notag
J_{1}\lesssim& g_{\tau}^{2}(s)
\int_{\mathbb{R}^{3}}\int_{B}
\frac{|\widehat{\psi}(\lambda_{1}+\lambda_{2})
\widehat{\psi}(\lambda_{3}+\lambda_{4})|^{2}}{|\lambda_{1}|^{1-\beta+2\delta}|\lambda_{2}\lambda_{3}\lambda_{4}|^{1-\beta}}
d\lambda_{1}d\lambda_{2}d\lambda_{3}d\lambda_{4}
\\\label{estimateJ1}\lesssim& g_{\tau}^{2}(s)
\int_{\mathbb{R}}
\frac{|\widehat{\psi}(u)|^{2}}{|u|^{1-2\beta+2\delta}}du
\int_{\mathbb{R}}
\frac{|\widehat{\psi}(v)|^{2}}{|v|^{1-2\beta}}
dv<\infty.
\end{align}
By Lemma \ref{lemma:L:ratio:diff}, for any $\delta>0$,
\begin{align}\notag
J_{2}\lesssim& g_{\tau}^{2}(s)
\int_{\mathbb{R}^{3}}\int_{B^{c}}
\frac{|\widehat{\psi}(\lambda_{1}+\lambda_{2})
\widehat{\psi}(\lambda_{3}+\lambda_{4})|^{2}}{|\lambda_{1}|^{1-\beta+4\tau-2\delta}|\lambda_{2}\lambda_{3}\lambda_{4}|^{1-\beta}}
d\lambda_{1}d\lambda_{2}d\lambda_{3}d\lambda_{4}
\\\label{estimateJ2}\lesssim& g_{\tau}^{2}(s)
\int_{\mathbb{R}}
\frac{|\widehat{\psi}(u)|^{2}}{|u|^{1-2\beta+4\tau-2\delta}}du
\int_{\mathbb{R}}
\frac{|\widehat{\psi}(v)|^{2}}{|v|^{1-2\beta}}
dv,
\end{align}
where the integrals are finite when $-\tau<(1-2\beta)/4$ and $\delta$ is sufficiently small.
By substituting the estimates (\ref{estimateJ1}) and (\ref{estimateJ2}) into (\ref{sumJ1J2}),
we obtain
\begin{align}\label{kol:distance3}
\mathbb{E}\left[|\widetilde{Z}^{(2)}_{s}(t)-Z^{(2)}(t)|^{2}\right]
\lesssim g_{\tau}^{2}(s)
\end{align}
under the condition $-\tau<(1-2\beta)/4$.
\\

\noindent{$\bullet$ \it Upper bound of $d_{\textup{Kol}}\left(Z(t)+\varepsilon,Z(t)\right)$}:
According to Definition \ref{def:Kol:distance} for the Kolmogorov distance,
\begin{align}\notag
d_{\textup{Kol}}\left(Z(t)+\varepsilon,Z(t)\right)
=&\underset{z\in \mathbb{R}}{\sup}\left|\mathbb{P}(Z(t)+\varepsilon\leq z)-\mathbb{P}(Z(t)\leq z)\right|
\\\notag=&\underset{z\in \mathbb{R}}{\sup}\left|\mathbb{P}(Z(t)\leq z-\varepsilon)-\mathbb{P}(Z(t)\leq z)\right|
\\\notag=&\underset{z\in \mathbb{R}}{\sup}\ \mathbb{P}(z-\varepsilon<Z(t)\leq z)
\\\label{levyineq0}=&\underset{z\in \mathbb{R}}{\sup}\ \mathbb{P}(z-\varepsilon/2<Z(t)\leq z+\varepsilon/2).
\end{align}
Because $Z(t)$ is a random variable defined as the sum of the second Wiener chaos, $Z^{(1)}(t)$, as defined in (\ref{def:Zlimit1}), and the fourth Wiener chaos, $Z^{(2)}(t)$, as defined in (\ref{def:Zlimit2}),
there exists a constant $C>0$ such that
\begin{align}\label{levyineq}
\mathbb{P}\left(\left|Z(t)-z\right|\leq \varepsilon/2 \right)\leq C \left\{\mathbb{E}\left[|Z(t)-z|^{2}\right]\right\}^{-1/8}\varepsilon^{1/4}
\end{align}
for any $\varepsilon>0$ and $z\in \mathbb{R}$.
The inequality (\ref{levyineq}) can be derived by using
the method employed in the proofs of Theorem 3.1 and Lemma 4.3 in \cite{MR3003367}, as well as the proof of Lemma 4 in \cite{MR4003569}.
For more details, please refer to \ref{sec:proof:levyineq}.
Substituting (\ref{levyineq}) into (\ref{levyineq0}) yields
\begin{align}\notag
d_{\textup{Kol}}\left(Z(t)+\varepsilon,Z(t)\right)
\leq& C \underset{z\in \mathbb{R}}{\sup} \left\{\mathbb{E}\left[|Z(t)-z|^{2}\right]\right\}^{-1/8}\varepsilon^{1/4}
\\\label{kol:distance4}=&C\left[\textup{Var}\left(Z(t)\right)\right]^{-1/8}\varepsilon^{1/4}.
\end{align}

\noindent{$\bullet$ \it Upper bound of $d_{\textup{Kol}}\left(Z(t)+Z^{(0)}_{s}-Z^{(0)},Z(t)\right)$}:
By (\ref{kol:distance4}),
\begin{align}\label{kol:distance5a}
d_{\textup{Kol}}\left(Z(t)+Z^{(0)}_{s}-Z^{(0)},Z(t)\right)
\lesssim |Z^{(0)}_{s}-Z^{(0)}|^{1/4}.
\end{align}
From (\ref{def:Z0s}) and Lemma \ref{lemma:chasodep:UY},
\begin{align}\notag
Z_{s}^{(0)}=&2^{4}s^{2\beta}L^{-2}(s)\int_{\mathbb{R}^{2}}|\widehat{\psi}(s(\lambda_{1}+\lambda_{2}))|^{2}
f_{X}(\lambda_{1})f_{X}(\lambda_{2})1_{\{\lambda_{1}\lambda_{2}<0\}}d\lambda_{1}d\lambda_{2}
\\\notag=&2^{4}\int_{\mathbb{R}^{2}}|\widehat{\psi}(\lambda_{1}+\lambda_{2})|^{2}
\frac{L(s/|\lambda_{1}|)}{L(s)}\frac{L(s/|\lambda_{2}|)}{L(s)}|\lambda_{1}\lambda_{2}|^{\beta-1}
1_{\{\lambda_{1}\lambda_{2}<0\}}d\lambda_{1}d\lambda_{2}.
\end{align}
By the definition of $Z^{(0)}$ in (\ref{def:Z0}),
\begin{align}\label{kol:distance5b}
|Z_{s}^{(0)}-Z^{(0)}|
\leq2^{4}\int_{\mathbb{R}^{2}}|\widehat{\psi}(\lambda_{1}+\lambda_{2})|^{2}
\left|\frac{L(s/|\lambda_{1}|)}{L(s)}\frac{L(s/|\lambda_{2}|)}{L(s)}
-1\right||\lambda_{1}\lambda_{2}|^{\beta-1}d\lambda_{1}d\lambda_{2}.
\end{align}
By the inequality
$$
\left|\frac{L(s/|\lambda_{1}|)}{L(s)}\frac{L(s/|\lambda_{2}|)}{L(s)}
-1\right|\leq \left|\frac{L^{2}(s/|\lambda_{1}|)}{L^{2}(s)}
-1\right|+\left|\frac{L^{2}(s/|\lambda_{2}|)}{L^{2}(s)}
-1\right|,
$$
for any $\delta>0$, we can estimate (\ref{kol:distance5b}) as follows
\begin{align}\notag
|Z_{s}^{(0)}-Z^{(0)}|
\leq&2^{5}\int_{\mathbb{R}^{2}}|\widehat{\psi}(\lambda_{1}+\lambda_{2})|^{2}
\left|\frac{L^{2}(s/|\lambda_{1}|)}{L^{2}(s)}
-1\right||\lambda_{1}\lambda_{2}|^{\beta-1}d\lambda_{1}d\lambda_{2}
\\\notag\lesssim&
g_{\tau}(s)\int_{\mathbb{R}}\int_{B}|\widehat{\psi}(\lambda_{1}+\lambda_{2})|^{2}
|\lambda_{1}|^{\beta-1-\delta}|\lambda_{2}|^{\beta-1}d\lambda_{1}d\lambda_{2}
\\\notag&+g_{\tau}(s)\int_{\mathbb{R}}\int_{B^{c}}|\widehat{\psi}(\lambda_{1}+\lambda_{2})|^{2}
|\lambda_{1}|^{\beta-1-2\tau+\delta}|\lambda_{2}|^{\beta-1}d\lambda_{1}d\lambda_{2},
\end{align}
where the last inequality follows from Lemma \ref{lemma:L:ratio:easy} and Lemma \ref{lemma:L:ratio:diff}.
According to Lemma \ref{lemma:integral:nu12}, when $2\beta<1$ and $-\tau<(1-2\beta)/2$, the integrals above are finite for sufficiently small $\delta>0$.
Hence,
$
|Z_{s}^{(0)}-Z^{(0)}| \lesssim g_{\tau}(s).
$
As a result,  (\ref{kol:distance5a}) implies that
\begin{align}\label{kol:distance5}
d_{\textup{Kol}}\left(Z(t)+Z^{(0)}_{s}-Z^{(0)},Z(t)\right)
\lesssim g_{\tau}^{1/4}(s).
\end{align}

Finally, by substituting the estimates (\ref{kol:distance2}), (\ref{kol:distance3}), (\ref{kol:distance4}), and (\ref{kol:distance5})
into (\ref{kol:distance1}),
we obtain
\begin{align}\label{terminal1}
d_{\textup{Kol}}\left(Z_{s}(t),Z(t)\right)
\lesssim& \varepsilon^{-2}g_{\tau}^{2}(s)+\varepsilon^{1/4}+g_{\tau}^{1/4}(s)
\end{align}
under the condition $-\tau<(1-2\beta)/4$.
By choosing $\varepsilon$ such that $\varepsilon^{-2}g_{\tau}^{2}(s)=\varepsilon^{1/4}$, achieved by setting
$\varepsilon=g^{8/9}_{\tau}(s)$, and noting that $g_{\tau}(s)\rightarrow0$ as $s\rightarrow\infty$,
we can rewrite (\ref{terminal1}) as
\begin{align*}
d_{\textup{Kol}}\left(Z_{s}(t),Z(t)\right)
\lesssim g^{2/9}_{\tau}(s)+g_{\tau}^{1/4}(s)\lesssim g^{2/9}_{\tau}(s).
\end{align*}


\subsection{Proof of Lemma \ref{lemma:integral:nu12}}\label{sec:proof:lemma:integral:nu12}

By H$\ddot{\textup{o}}$lder's inequality,
\begin{align*}
\left|\int_{\mathbb{R}}
\frac{|\widehat{\psi}(\lambda_{1}+\xi)\widehat{\psi}(\lambda_{2}-\xi)|}{\left|\xi\right|^{1-\beta-\nu_{3}}} d\xi
\right|^{2}
\leq
\left(\int_{\mathbb{R}}
\frac{|\widehat{\psi}(\lambda_{1}+\xi)|^{2}}{\left|\xi\right|^{1-\beta-\nu_{3}}} d\xi \right)
\left(\int_{\mathbb{R}}
\frac{|\widehat{\psi}(\lambda_{2}-\xi)|^{2}}{\left|\xi\right|^{1-\beta-\nu_{3}}} d\xi\right).
\end{align*}
Hence,
\begin{align*}
&\int_{\mathbb{R}^{2}}
\left|\lambda_{1}\right|^{\nu_{1}+\beta-1}\left|\lambda_{2}\right|^{\nu_{2}+\beta-1}
\left|\int_{\mathbb{R}}
\frac{|\widehat{\psi}(\lambda_{1}+\xi)\widehat{\psi}(\lambda_{2}-\xi)|}{\left|\xi\right|^{1-\beta-\nu_{3}}} d\xi
\right|^{2}
d\lambda_{1}d\lambda_{2}
\\\leq&
\left(\int_{\mathbb{R}}
\left|\lambda_{1}\right|^{\nu_{1}+\beta-1}
\int_{\mathbb{R}}
\frac{|\widehat{\psi}(\lambda_{1}+\xi)|^{2}}{\left|\xi\right|^{1-\beta-\nu_{3}}} d\xi
d\lambda_{1}\right)
\left(\int_{\mathbb{R}}
\left|\lambda_{2}\right|^{\nu_{2}+\beta-1}
\int_{\mathbb{R}}
\frac{|\widehat{\psi}(\lambda_{2}-\xi)|^{2}}{\left|\xi\right|^{1-\beta-\nu_{3}}} d\xi
d\lambda_{2}\right)
\\=&\left(\int_{\mathbb{R}}
\frac{|\widehat{\psi}(z)|^{2}}{\left|z\right|^{1-2\beta-(\nu_{1}+\nu_{3})}} dz
\right)\left(\int_{\mathbb{R}}
\frac{|\widehat{\psi}(z)|^{2}}{\left|z\right|^{1-2\beta-(\nu_{2}+\nu_{3})}} dz
\right)<\infty,
\end{align*}
where the last equal sign follows from the Riesz composition formula \cite[Appendix A]{anh2003higher},
which holds under the conditions $0<\nu_{k}+\beta<1$ for $k\in\{1,2,3\}$ and  $0<\nu_{\ell}+\nu_{3}+2\beta<1$ for $\ell\in\{1,2\}$.

\subsection{Proof of Theorem \ref{thm:Stein}}\label{sec:proof:thm:Stein}
By the definition of $U^{2}[s]Y$ in (\ref{def:;scalogram}),
$
sU^{2}[s]Y(st_{k}) = (F_{2k-1,s})^{2}+(F_{2k,s})^{2},\ k=1,2,...,n,
$
where $F_{2k-1,s} = s^{1/2}Y\star \psi_{s}(st_{k})$ and $F_{2k,s} = s^{1/2}Y\star \mathcal{H}\psi_{s}(st_{k})$.
By Lemma \ref{lemma:chisquare} and the identity $\widehat{\mathcal{H}\psi}(\lambda) = -i\textup{sgn}(\lambda)\widehat{\psi}(\lambda)$,
we can express $F_{2k-1,s}$ and $F_{2k,s}$ as follows:
\begin{align}\notag
F_{\bullet,s} = \int_{\mathbb{R}^{2}}^{'}f_{\bullet,s}(\lambda_{1},\lambda_{2})W(d\lambda_{1})W(d\lambda_{2}),
\end{align}
where
\begin{align}\label{def:f2k-1}
f_{2k-1,s}(\lambda_{1},\lambda_{2}) = 2s^{1/2}
e^{ist_{k}(\lambda_{1}+\lambda_{2})} \widehat{\psi}(s(\lambda_{1}+\lambda_{2})) \sqrt{f_{X}(\lambda_{1})f_{X}(\lambda_{2})}1_{\{\lambda_{1}\lambda_{2}<0\}}
\end{align}
and
\begin{align}\notag
f_{2k,s}(\lambda_{1},\lambda_{2}) = -i2s^{1/2}
e^{ist_{k}(\lambda_{1}+\lambda_{2}) }\textup{sgn}(\lambda_{1}+\lambda_{2})\widehat{\psi}(s(\lambda_{1}+\lambda_{2})) \sqrt{f_{X}(\lambda_{1})f_{X}(\lambda_{2})}1_{\{\lambda_{1}\lambda_{2}<0\}}.
\end{align}
Let $\mathbf{F}_{s} = (F_{1,s},F_{2,s},...,F_{2n-1,s},F_{2n,s})\in \mathbb{R}^{2n}$, and define $\mathbf{C}_{s}\in \mathbb{R}^{2n\times 2n}$ as its covariance matrix.
Additionally, let $\mathbf{N}_{s}$ be a $2n$-dimensional normal random vector
with mean zero and covariance matrix $\mathbf{C}_{s}= [C_{s}(\ell,\ell')]_{1\leq \ell,\ell'\leq 2n}$, represented as $(N_{1,s},N_{2,s},...,N_{2n-1,s},N_{2n,s})$.
Recalling the definition of the Wasserstein distance between $\mathbf{F}_{s}$
and $\mathbf{N}_{s}$:
\begin{align*}
d_{\textup{W}}(\mathbf{F}_{s},\mathbf{N}_{s})
= \underset{\|h\|_{\textup{Lip}}\leq 1}{\textup{sup}}\left|\mathbb{E}[h(\mathbf{F}_{s})]-\mathbb{E}[h(\mathbf{N}_{s})]\right|,
\end{align*}
where the supremum is taken over continuous functions $h: \mathbb{R}^{2n}\rightarrow \mathbb{R}$ with
$|h(\mathbf{x})-h(\mathbf{y})|\leq |\mathbf{x}-\mathbf{y}|$ for $\mathbf{x},\mathbf{y}\in \mathbb{R}^{2n}$.
For any function  $g:\mathbb{R}\rightarrow\mathbb{R}$ satisfying $\|g\|_{\textup{Lip}}\leq 1$ and any $(c_{1},...,c_{n})\in \mathbb{R}^{n}$ with $c_{1}^{2}+\cdots +c_{n}^{2}\leq 1$,
the mapping
$$
(x_1,x_2,...,x_{2n-1},x_{2n})\in \mathbb{R}^{2n}\mapsto
g\left(\overset{n}{\underset{k=1}{\sum}}c_{k}\sqrt{x_{2k-1}^{2}+x_{2k}^{2}}\right)
$$
is continuous with a Lipschitz constant less than one.
Hence,
\begin{align}\notag
&d_{\textup{W}}\left(\overset{n}{\underset{k=1}{\sum}}c_{k}s^{1/2}U[s]Y(st_{k}),
\overset{n}{\underset{k=1}{\sum}}c_{k}\sqrt{N_{2k-1,s}^{2}+N_{2k,s}^{2}}\right)
\\\notag=&d_{\textup{W}}\left(\overset{n}{\underset{k=1}{\sum}}c_{k}\sqrt{F_{2k-1,s}^{2}+F_{2k,s}^{2}},
\overset{n}{\underset{k=1}{\sum}}c_{k}\sqrt{N_{2k-1,s}^{2}+N_{2k,s}^{2}}\right)
\\\notag=& \underset{\|g\|_{\textup{Lip}}\leq 1}{\textup{sup}}\left|\mathbb{E}\left[g\left(\overset{n}{\underset{k=1}{\sum}}c_{k}\sqrt{F_{2k-1,s}^{2}+F_{2k,s}^{2}}\right)\right]-
\mathbb{E}\left[g\left(\overset{n}{\underset{k=1}{\sum}}c_{k}\sqrt{N_{2k-1,s}^{2}+N_{2k,s}^{2}}\right)\right]\right|
\\\label{def:dW2}\leq&
d_{\textup{W}}(\mathbf{F}_{s},\mathbf{N}_{s}).
\end{align}
According to \cite[Theorem 6.1.1]{nourdin2012normal},
\begin{align}\label{multivariate:Stein}
d_{\textup{W}}(\mathbf{F}_{s},\mathbf{N}_{s})
\leq \sqrt{2}
\|\mathbf{C}_{s}^{-1}\|_{\textup{op}}
\|\mathbf{C}_{s}\|_{\textup{op}}^{1/2}
\sqrt{\overset{2n}{\underset{\ell,\ell'=1}{\sum}}
\mathbb{E}\left[\left(C_{s}(\ell,\ell')-\langle DF_{\ell',s},-DL^{-1}F_{\ell,s}\rangle_{\overline{H}}\right)^{2}\right]},
\end{align}
where $\|\mathbf{C}_{s}^{-1}\|_{\textup{op}}$ represents the operator norm of the inverse matrix of $\mathbf{C}_{s}$, $D$ denotes the Malliavin derivative, $L^{-1}$ signifies the pseudo-inverse of the infinitesimal generator $L$ for the Ornstein-Uhlenbeck semigroup,
and $\overline{H}=\{f\in L^{2}(\mathbb{R})\mid f(-\lambda)=\overline{f(\lambda)}\ \textup{for all}\ \lambda\in \mathbb{R}\}$ is a complex Hilbert space equipped
with the inner product
$
\langle f,g\rangle_{\overline{H}} = \int_{\mathbb{R}}f(\lambda)\overline{g(\lambda)}d\lambda.
$

Denote
\begin{align*}
f^{*2}(u) = \int_{\mathbb{R}}f_{X}(u-v)f_{X}(v)1_{\{(u-v)v<0\}}dv,\ u\in \mathbb{R},
\end{align*}
which is an even and continuous function.
By the orthogonal property (\ref{ortho}), for any $\ell,\ell'\in\{1,2,...,n-1,2n\}$,
\begin{equation}\label{def:Cs:structure}
\frac{C_{s}(\ell,\ell')}{16} = \left\{\begin{array}{l}
\int_{0}^{\infty}\cos((t_{k}-t_{k'})u)|\widehat{\psi}(u)|^{2}f^{*2}(u/s)du\ \textup{if}\ (\ell,\ell')=(2k-1,2k'-1)\ \textup{or}\ (2k,2k'),
\\
\int_{0}^{\infty}\sin((t_{k}-t_{k'})u)|\widehat{\psi}(u)|^{2}f^{*2}(u/s)du\ \textup{if}\ (\ell,\ell')=(2k-1,2k'),
\\
\int_{0}^{\infty}\sin((t_{k'}-t_{k})u)|\widehat{\psi}(u)|^{2}f^{*2}(u/s)du\ \textup{if}\ (\ell,\ell')=(2k,2k'-1)
\end{array}
\right.
\end{equation}
for $k\in\{1,2,...,n-1,,n\}$.
Define a matrix $\mathbf{C}^{(n)}_{s}\in \mathbb{C}^{n\times n}$, whose components are given by
\begin{equation}\label{def:Phi_complex}
C^{(n)}_{s}(k,k') =
16\int_{0}^{\infty}e^{i(t_{k}-t_{k'})u}|\widehat{\psi}(u)|^{2}f^{*2}(u/s)du,\ k,k'=1,2,...,n.
\end{equation}
For any $\mathbf{v} = [v_{1}\ v_{2}\ \cdots\ v_{2n-1}\ v_{2n}]\in \mathbb{R}^{2n},$
we denote
$\mathbf{v}^{(n)} = [v_{1}+iv_{2}\ \cdots\ v_{2n-1}+iv_{2n}]\in \mathbb{C}^{n}.$
We notice that the structure (\ref{def:Cs:structure}) for the matrix $\mathbf{C}_{s}$ enables us to derive the following equality
\begin{align}\label{real_complex_link}
\mathbf{v}\mathbf{C}_{s} \mathbf{v}^{\top}= \mathbf{v}^{(n)}\mathbf{C}^{(n)}_{s}(\mathbf{v}^{(n)})^{*},
\end{align}
where $\top$ denotes the transpose and $*$ represents the complex conjugate.
Because there exists a threshold $\widetilde{s}$ such that  $|\widehat{\psi}(\cdot)|^{2}f^{*2}(\cdot/s)1_{[0,\infty)}(\cdot)$ is a nonvanishing, nonnegative and  integrable function for any $s>\widetilde{s}$,
the Fourier transform of $|\widehat{\psi}(\cdot)|^{2}f^{*2}(\cdot/s)1_{[0,\infty)}(\cdot)$ is strictly positive definite \cite[Corollary 6.9]{wendland2004scattered} at least for $s>\widetilde{s}$.
Hence, $\mathbf{C}^{(n)}_{s}$ is positive definite for $s>\widetilde{s}$.
This implies that for any nonzero vector  $\mathbf{v}^{(n)}\in \mathbb{C}^{n}$,
$ \mathbf{v}^{(n)}\mathbf{C}^{(n)}_{s}(\mathbf{v}^{(n)})^{*}>0$ for $s>\widetilde{s}$.
This property, combined with the observation (\ref{real_complex_link}), implies the positive definiteness of $\mathbf{C}_{s}$ for $s>\widetilde{s}$,
thus ensuring that $\|\mathbf{C}_{s}^{-1}\|_{\textup{op}}$ in (\ref{multivariate:Stein}) is finite for $s>\widetilde{s}$.

Under the condition $2\beta>1$, for any $u\in \mathbb{R}$, $f^{*2}(u/s)\rightarrow \|f_{X}\|_{2}^{2}$ as $s\rightarrow\infty$.
Hence,
$$\underset{s\rightarrow\infty}{\lim}
C_{s}(\ell,\ell') = C_{\infty}(\ell,\ell'):=16\|f_{X}\|_{2}^{2}\Psi(\ell,\ell'),\ \ell,\ell'\in\{1,2,...,2n-1,2n\},
$$
where $\Psi$ is defined in (\ref{def:Phi:bigmatrix}). The positive definiteness of $\Psi$ can also be proved using the observation (\ref{real_complex_link}),
along with the fact $0<\int_{0}^{\infty} |\widehat{\psi}(u)|^{2}du<\infty$ and the strictly positive definiteness
of the Fourier transform of  $|\widehat{\psi}(\cdot)|^{2}$.
Hence, $\|\mathbf{C}_{s}^{-1}\|_{\textup{op}}
\|\mathbf{C}_{s}\|_{\textup{op}}^{1/2}$ converges when $s\rightarrow\infty$ and (\ref{multivariate:Stein}) implies that
\begin{align}\label{multivariate:Stein2}
d_{\textup{W}}(\mathbf{F}_{s},\mathbf{N}_{s})
\lesssim\sqrt{\overset{2n}{\underset{\ell,\ell'=1}{\sum}}
\mathbb{E}\left[\left(C_{s}(\ell,\ell')-\langle DF_{\ell,s},-DL^{-1}F_{\ell',s}\rangle_{\overline{H}}\right)^{2}\right]}.
\end{align}

Next, we aim to estimate the right hand side of (\ref{multivariate:Stein2}). For any $p\in \mathbb{N}$ and any function $f$ in the $p$-th tensor product of $\overline{H}$, denoted as $\overline{H}^{\otimes p}$, let
$I_{p}(f)$ represent the $p$-fold Wiener integral of $f$ with respect to the complex-valued Gaussian white noise random measure $W$.
According to \cite[Chapter 2]{nourdin2012normal}, for $\ell,\ell'\in\{1,...,2n\}$,
\begin{align}\label{lemma:DIp}
DF_{\ell,s} = DI_{2}(f_{\ell,s})
= 2 I_{1}(f_{\ell,s}),
\end{align}
where
\begin{align*}
I_{1}(f_{\ell,s})(\cdot) = \int_{\mathbb{R}}f_{\ell,s}(\lambda,\cdot)
W(d\lambda)
\end{align*}
and
\begin{align}\label{lemma:L-1Ip}
L^{-1}I_{2}(f_{\ell',s}) = -\frac{1}{2} I_{2}(f_{\ell',s}).
\end{align}
By (\ref{lemma:DIp}), (\ref{lemma:L-1Ip}), and the product formula \cite{major1981lecture},
\begin{align}\notag
\langle DF_{\ell,s},-DL^{-1}F_{\ell',s}\rangle_{\overline{H}}
=&\langle 2I_{1}(f_{\ell,s}),I_{1}(f_{\ell',s})\rangle_{\overline{H}}
\\\label{compute:innerproduct}=&2\langle f_{\ell,s},f_{\ell',s}\rangle_{\overline{H}^{\otimes 2}}+2I_{2}(f_{\ell,s}\widetilde{\otimes}_{1} f_{\ell',s}),
\end{align}
where
$$
\langle f_{\ell,s},f_{\ell',s}\rangle_{\overline{H}^{\otimes 2}}=\int_{\mathbb{R}^{2}} f_{\ell,s}(\lambda_{1},\lambda_{2})\overline{f_{\ell',s}(\lambda_{1},\lambda_{2})}\,d\lambda_{1}d\lambda_{2},
$$
\begin{align}\label{def:contraction}
f_{\ell,s}\otimes_{1} f_{\ell',s}(\lambda_{1},\lambda_{2})
= \int_{\mathbb{R}} f_{\ell,s}(\lambda_{1},\eta)
f_{\ell',s}(\lambda_{2},-\eta)\, d\eta,
\end{align}
and $f_{\ell,s}\widetilde{\otimes}_{1} f_{\ell',s}$ is the canonical symmetrization of the product contraction $f_{\ell,s}\otimes_{1} f_{\ell',s}$, defined as
\begin{align*}
f_{\ell,s}\widetilde{\otimes}_{1} f_{\ell',s}(\lambda_{1},\lambda_{2})
=2^{-1}\left[ f_{\ell,s}\otimes_{1} f_{\ell',s}(\lambda_{1},\lambda_{2})+f_{\ell,s}\otimes_{1} f_{\ell',s}(\lambda_{2},\lambda_{1})\right].
\end{align*}
On the other hand, for $\ell,\ell'\in\{1,...,2n\}$, by the isometry property of the multiple Wiener-It$\hat{\textup{o}}$ integrals,
\begin{align}\label{compute:Cs}
C_{s}(\ell,\ell') = \mathbb{E}\left[I_{2}(f_{\ell,s})I_{2}(f_{\ell',s})\right] = 2\langle f_{\ell,s},f_{\ell',s}\rangle_{\overline{H}^{\otimes 2}}.
\end{align}
Combining (\ref{compute:innerproduct}) and (\ref{compute:Cs}) yields
\begin{align}\notag
&\mathbb{E}\left[\left(C_{s}(\ell,\ell')-\langle DF_{\ell,s},-DL^{-1}F_{\ell',s}\rangle_{\overline{H}}\right)^{2}\right]
\\\notag=&\mathbb{E}\left[\left(2I_{2}(f_{\ell,s}\widetilde{\otimes}_{1} f_{\ell',s})\right)^{2}\right]
\\\notag=&8 \langle f_{\ell,s}\widetilde{\otimes}_{1} f_{\ell',s},f_{\ell,s}\widetilde{\otimes}_{1} f_{\ell',s}\rangle_{\overline{H}^{\otimes 2}}
\\\label{estimate:Cinner1}\leq& 8\langle f_{\ell,s}\otimes_{1} f_{\ell',s},f_{\ell,s}\otimes_{1} f_{\ell',s}\rangle_{\overline{H}^{\otimes 2}},
\end{align}
where the inequality follows from the Minkowski inequality.
Denote
\begin{align*}
f_{s}(\lambda_{1},\lambda_{2}) =
s^{1/2}|\widehat{\psi}(s(\lambda_{1}+\lambda_{2}))| \sqrt{f_{X}(\lambda_{1})f_{X}(\lambda_{2})}.
\end{align*}
According to the definition of $f_{\ell,s}$ in (\ref{def:f2k-1}), $f_{s}\geq |f_{\ell,s}|/2$ for $\ell\in\{1,...,2n\}$.
Hence, (\ref{estimate:Cinner1}) implies that
\begin{align}\label{estimate:Cinner2}
\overset{2n}{\underset{\ell,\ell'=1}{\sum}}\mathbb{E}\left[\left(C_{s}(\ell,\ell')-\langle DF_{\ell,s},-DL^{-1}F_{\ell',s}\rangle_{\overline{H}}\right)^{2}\right]
\leq 2^{9}n^{2} \|f_{s}\otimes_{1}f_{s}\|_{\overline{H}^{\otimes 2}}^{2}.
\end{align}
Using (\ref{def:dW2}) and substituting (\ref{estimate:Cinner2}) into (\ref{multivariate:Stein2}), we obtain
 \begin{align*}
d_{\textup{W}}\left(\overset{n}{\underset{k=1}{\sum}}c_{k}s^{1/2}U[s]Y(st_{k}),
\overset{n}{\underset{k=1}{\sum}}c_{k}\sqrt{N_{2k-1,s}^{2}+N_{2k,s}^{2}}\right)
\lesssim \|f_{s}\otimes_{1}f_{s}\|_{\overline{H}^{\otimes 2}}.
\end{align*}
To complete the proof, we need to show that
$\|f_{s}\otimes_{1}f_{s}\|_{\overline{H}^{\otimes 2}}\lesssim s^{-1/2}$ for the case $\beta\in(3/4,1)$
and
$\|f_{s}\otimes_{1}f_{s}\|_{\overline{H}^{\otimes 2}}\lesssim s^{-1/2+2(3/4-\beta)}L^{2}(s)$ for the case $\beta\in(1/2,3/4)$.
From the definition of the product contraction $\otimes_{1}$ in (\ref{def:contraction}),
\begin{align}\notag
\|f_{s}\otimes_{1}f_{s}\|_{\overline{H}^{\otimes 2}}^{2}
=&s^{2}\int_{\mathbb{R}^{2}}
f_{X}(\lambda_{1})f_{X}(\lambda_{2})
\left[\int_{\mathbb{R}}|\widehat{\psi}(s(\lambda_{1}+\eta_{1}))| |\widehat{\psi}(s(\lambda_{2}-\eta_{1}))|
f_{X}(\eta_{1})d\eta_{1}\right]
\\\notag&\times\left[\int_{\mathbb{R}}|\widehat{\psi}(s(\lambda_{1}+\eta_{2}))| |\widehat{\psi}(s(\lambda_{2}-\eta_{2}))|
f_{X}(\eta_{2})d\eta_{2}\right]d\lambda_{1}d\lambda_{2}
\\\notag=&
s^{-1}\int_{\mathbb{R}^{3}}
\left[\int_{\mathbb{R}}f_{X}(\lambda)f_{X}(\frac{x+y}{s}-\lambda) f_{X}(\frac{x}{s}-\lambda)f_{X}(\frac{z}{s}-\lambda)d\lambda\right]
\\\label{estimate:singularintegral1}&\times|\widehat{\psi}(x)| |\widehat{\psi}(y)||\widehat{\psi}(z)| |\widehat{\psi}(x+y-z)|
\,dx\,dy\, dz.
\end{align}

\noindent{$\bullet$} For the case $\beta\in (3/4,1)$ (i.e., $4(\beta-1)>-1$), because $f_{X}(\lambda) = L(|\lambda|^{-1})|\lambda|^{\beta-1}$, and for any $\delta>0$, there exists a constant $C>0$ such that
$L(|\lambda|^{-1})\leq C |\lambda|^{-\delta}$ for all $\lambda\in(-1,1)\setminus\{0\}$,  it follows that $\|f_{X}\|_{4}<\infty$.
By the generalized H$\ddot{\textup{o}}$lder inequality,
\begin{align}\label{generalHolder}
\int_{\mathbb{R}}f_{X}(\lambda)f_{X}(\frac{x+y}{s}-\lambda) f_{X}(\frac{x}{s}-\lambda)f_{X}(\frac{z}{s}-\lambda)d\lambda
\leq \|f_{X}\|_{4}^{4}.
\end{align}
By applying (\ref{generalHolder}) to (\ref{estimate:singularintegral1}), we obtain
\begin{align*}
\|f_{s}\otimes_{1}f_{s}\|_{\overline{H}^{\otimes 2}}^{2}
\lesssim s^{-1}.
\end{align*}

\noindent{$\bullet$} For the case $\beta\in(1/2,3/4)$,
by considering a change of variable $\lambda = u/s$, the expression inside the brackets in (\ref{estimate:singularintegral1}) can be rewritten as
\begin{align}\notag
&\int_{\mathbb{R}}f_{X}(\lambda)f_{X}(\frac{x+y}{s}-\lambda) f_{X}(\frac{x}{s}-\lambda)f_{X}(\frac{z}{s}-\lambda)d\lambda
\\\notag=&
s^{-1}\int_{\mathbb{R}}f_{X}(\frac{u}{s})f_{X}(\frac{x+y-u}{s}) f_{X}(\frac{x-u}{s})f_{X}(\frac{z-u}{s})du
\\\notag=&
s^{-1+4(1-\beta)}\int_{\mathbb{R}}\frac{L(s|u|^{-1})}{|u|^{1-\beta}}
\frac{L(s|x+y-u|^{-1})}{|x+y-u|^{1-\beta}}
\frac{L(s|x-u|^{-1})}{|x-u|^{1-\beta}}
\frac{L(s|z-u|^{-1})}{|z-u|^{1-\beta}}
du
\\\notag=&
s^{-1+4(1-\beta)}L^{4}(s)\int_{\mathbb{R}}
\left[\frac{L(s|u|^{-1})}{L(s)}
\frac{L(s|x+y-u|^{-1})}{L(s)}
\frac{L(s|x-u|^{-1})}{L(s)}
\frac{L(s|z-u|^{-1})}{L(s)}\right]
\\\label{estimate:singularintegral3}&\times \frac{1}{
|u(x+y-u)(x-u)(z-u)|^{1-\beta}}
du.
\end{align}
By substituting (\ref{estimate:singularintegral3}) into (\ref{estimate:singularintegral1}),
we obtain
\begin{align}\label{estimate:fsotimesfs}
\|f_{s}\otimes_{1}f_{s}\|_{\overline{H}^{\otimes 2}}^{2}
=s^{-2+4(1-\beta)}L^{4}(s) \mathcal{K}_{s},
\end{align}
where
\begin{align}\notag
\mathcal{K}_{s}=&
\int_{\mathbb{R}^{4}}
\left[\frac{L(s|u|^{-1})}{L(s)}
\frac{L(s|x+y-u|^{-1})}{L(s)}
\frac{L(s|x-u|^{-1})}{L(s)}
\frac{L(s|z-u|^{-1})}{L(s)}\right]
\\\notag&\times \frac{|\widehat{\psi}(x)| |\widehat{\psi}(y)||\widehat{\psi}(z)| |\widehat{\psi}(x+y-z)|}{
|u(x+y-u)(x-u)(z-u)|^{1-\beta}}
\,du\,dx\,dy\, dz.
\end{align}
Let
$\mu=(\mu_{1},\mu_{2},\mu_{3},\mu_{4})\in\{-1,1\}^{4}$
be a function of $(|u|,|x+y-u|,|x-u|,|z-u|)$  defined as follows
\begin{align}\notag
&\begin{aligned}
&\left\{\begin{array}{ll}
\mu_{1}=1\  &\textup{if}\ |u|>1,
\\
\mu_{1}=-1\  &\textup{if}\ |u|\leq1,
\end{array}\right.
\\
&\left\{\begin{array}{ll}
\mu_{2}=1\  &\textup{if}\ |x+y-z|>1,
\\
\mu_{2}=-1\  &\textup{if}\ |x+y-u|\leq1,
\end{array}\right.
\end{aligned}
&\begin{aligned}
&\left\{\begin{array}{ll}
\mu_{3}=1\  &\textup{if}\ |x-u|>1,
\\
\mu_{3}=-1\  &\textup{if}\ |x-u|\leq1,
\end{array}\right.
\\
&\left\{\begin{array}{ll}
\mu_{4}=1\  &\textup{if}\ |z-u|>1,
\\
\mu_{4}=-1\  &\textup{if}\ |z-u|\leq1.
\end{array}\right.
\end{aligned}
\end{align}
By denoting $$R_{\mu}=\{(u,x,y,z)\in \mathbb{R}^{4}\mid  |u|^{\mu_{1}}\geq 1,\ |x+y-u|^{\mu_{2}}\geq 1,\ |x-u|^{\mu_{3}}\geq 1,\ |z-u|^{\mu_{4}}\geq 1\},$$
we can express $\mathcal{K}_{s}$  as
\begin{align*}
\mathcal{K}_{s} = \underset{\mu\in\{-1,1\}^{4}}{\sum}\mathcal{K}_{s,\mu},
\end{align*}
where
\begin{align}\notag
\mathcal{K}_{s,\mu}=&
\int_{\mathbb{R}_{\mu}}
\left[\frac{L(s|u|^{-1})}{L(s)}
\frac{L(s|x+y-u|^{-1})}{L(s)}
\frac{L(s|x-u|^{-1})}{L(s)}
\frac{L(s|z-u|^{-1})}{L(s)}\right]
\\\label{estimate:singularintegral4}&\times \frac{|\widehat{\psi}(x)| |\widehat{\psi}(y)||\widehat{\psi}(z)| |\widehat{\psi}(x+y-z)|}{
|u(x+y-u)(x-u)(z-u)|^{1-\beta}}
\,du\,dx\,dy\, dz.
\end{align}
According to \cite[Theorem 1.5.3]{bingham1989regular} (see also \cite[p. 1474]{leonenko2014sojourn}), for any $\delta>0$, there exist a constant $C>0$ and a threshold $s^{*}>0$
such that for any $s$ greater than $s^*$, the following inequalities hold
\begin{align}\label{farfield}
\underset{0<r<1}{\sup} \frac{(rs)^{\delta}L\left(rs\right)}{s^{\delta}L(s)}
\leq C
\end{align}
and
\begin{align}\label{nearfield}
\underset{r\geq1}{\sup} \frac{(rs)^{-\delta}L\left(rs\right)}{s^{-\delta}L(s)}
\leq C.
\end{align}
By applying (\ref{farfield}) and (\ref{nearfield}) to (\ref{estimate:singularintegral4}), we obtain the estimate
\begin{align*}
\mathcal{K}_{s,\mu}\leq
C^{4}\int_{\mathbb{R}_{\mu}}
 \frac{|\widehat{\psi}(x)| |\widehat{\psi}(y)||\widehat{\psi}(z)| |\widehat{\psi}(x+y-z)|}{
|u|^{1-\beta-\delta \mu_{1}}|x+y-u|^{1-\beta-\delta\mu_{2}}|x-u|^{1-\beta-\delta\mu_{3}}|z-u|^{1-\beta-\delta\mu_{4}}}
\,du\,dx\,dy\, dz
\end{align*}
for any $s$ greater than $s^{*}$. By Lemma \ref{lemma:singular_integral}, the integral above is finite for sufficiently small $\delta>0$.  
Therefore, from (\ref{estimate:fsotimesfs}), we obtain
\begin{align*}
\|f_{s}\otimes_{1}f_{s}\|_{\overline{H}^{\otimes 2}}^{2}
\lesssim s^{-2+4(1-\beta)}L^{4}(s)=s^{-1+4(3/4-\beta)}L^{4}(s)
\end{align*}
for the case $\beta\in(1/2,3/4)$.

\subsection{Proof of Lemma \ref{lemma:singular_integral}}\label{sec:proof:lemma:singular_integral}

For each $k\in\{1,2,3,4\}$, we denote $\beta_{k} = \beta+\delta\mu_{k}$.
Because $\beta\in(1/2,3/4)$, there exists a sufficiently small $\delta>0$ such that
$\beta_{k}\in(1/2,3/4)$ for all $k\in\{1,2,3,4\}$.
Let $\beta_{k}^{+}$ be a constant greater than $\beta_{k}$ and sufficiently close to $\beta_{k}$.
 By H$\ddot{\textup{o}}$lder's inequality,
\begin{align*}
\int_{\mathbb{R}}
 \frac{|\widehat{\psi}(z)| |\widehat{\psi}(x+y-z)|}{
|z-u|^{1-\beta_{4}}}
dz
=I^{1/2}_{\beta_{4}}|\widehat{\psi}|^{2}(u)
\times I^{1/2}_{\beta_{4}}|\widehat{\psi}|^{2}(x+y-u),
\end{align*}
where
\begin{align*}
I_{\beta_{4}}|\widehat{\psi}|^{2}(w):=
\int_{\mathbb{R}}
 \frac{|\widehat{\psi}(z)|^{2}}{
|w-z|^{1-\beta_{4}}}
dz,
\end{align*}
which is an even function because $\widehat{\psi}(-z)=\overline{\widehat{\psi}(z)}.$
Hence, (\ref{def:singular_integral}) can be bounded as follows
\begin{align}\label{estimate:singularintegral7}
\mathcal{K}\leq
\int_{\mathbb{R}^{3}}
 \frac{|\widehat{\psi}(x)| |\widehat{\psi}(y)|}{
|u|^{1-\beta_{1}}|x+y-u|^{1-\beta_{2}}|x-u|^{1-\beta_{3}}}\,
 I^{1/2}_{\beta_{4}}|\widehat{\psi}|^{2}(u)\,
I^{1/2}_{\beta_{4}}|\widehat{\psi}|^{2}(x+y-u)
\,du\,dx\,dy.
\end{align}
By denoting
\begin{align}\label{def:Q}
Q(w)=\int_{\mathbb{R}}|\widehat{\psi}(y)|
 \frac{I^{1/2}_{\beta_{4}}|\widehat{\psi}|^{2}(w+y)}{
|w+y|^{1-\beta_{2}}}
\,dy,\ w\in\mathbb{R},
\end{align}
we can rewrite (\ref{estimate:singularintegral7}) as
\begin{align}\label{estimate:singularintegral8}
\mathcal{K}\leq&
\int_{\mathbb{R}}\left[\int_{\mathbb{R}}
  \frac{ I^{1/2}_{\beta_{4}}|\widehat{\psi}|^{2}(u)}
{|u|^{1-\beta_{1}}}
\frac{Q(x-u)}{|x-u|^{1-\beta_{3}}}
\,du
\right]
 |\widehat{\psi}(x)|\,dx.
\end{align}
We will show that
\begin{equation}\tag{Claim 1}
 \frac{I^{1/2}_{\beta_{4}}|\widehat{\psi}|^{2}(\bullet)}{|\bullet|^{1-\beta_{1}}}\in L^{2/(3-2\beta_{1}^{+}-\beta_{4}^{+})},\label{Claim 1}
\end{equation}
\begin{equation}\tag{Claim 2}
Q\in  L^{2/(3-2\beta_{2}^{+}-\beta_{4}^{+})} \cap C(\mathbb{R}), \label{Claim 2}
\end{equation}
where $C(\mathbb{R})$ is the set of continuous functions, and
\begin{equation}\tag{Claim 3}
 \frac{Q(\bullet)}{|\bullet|^{1-\beta_{3}}}\in L^{2/(5-2\beta_{2}^{+}-2\beta_{3}^{+}-\beta_{4}^{+})}.\label{Claim 3}
\end{equation}

Denote $p = 2/(3-2\beta_{1}^{+}-\beta_{4}^{+})$ and $r=2/(5-2\beta_{2}^{+}-2\beta_{3}^{+}-\beta_{4}^{+})$.
Because $\beta\in (1/2,3/4)$, we have $p>4/3$ and $r>4/5$.
We consider the following two cases according to the range of $r$.

\noindent{\it Case 1:} For $r\in (4/5,1]$, \ref{Claim 2} and \ref{Claim 3} imply that
$Q(\bullet)|\bullet|^{\beta_{3}-1}\in L^{1}$.
By Minkowski's inequality for integrals, we have the following estimate for the convolution inside the bracket in (\ref{estimate:singularintegral8})
\begin{align}\label{result:case1}
\left\|\int_{\mathbb{R}}
  \frac{ I^{1/2}_{\beta_{4}}|\widehat{\psi}|^{2}(u)}
{|u|^{1-\beta_{1}}}
\frac{Q(\bullet-u)}{|\bullet-u|^{1-\beta_{3}}}
\,du\right\|_{p}\leq
 \left\|\frac{I^{1/2}_{\beta_{4}}|\widehat{\psi}|^{2}(\bullet)}{|\bullet|^{1-\beta_{1}}}\right\|_{p}
 \left\|\frac{Q(\bullet)}{|\bullet|^{1-\beta_{3}}}\right\|_{1}.
\end{align}

\noindent{\it Case 2:} For $r>1$,
let $q$ be a constant defined by
\begin{equation*}
\frac{1}{q} = \frac{1}{p}+\frac{1}{r}-1=3-(\beta_{1}^{+}+\beta_{2}^{+}+\beta_{3}^{+}+\beta_{4}^{+}).
\end{equation*}
Because $\beta\in(1/2,3/4)$, when $\delta$ is sufficiently small
and $\beta_{k}^{+}$ is sufficiently close to $\beta_{k}$ for each $k\in\{1,2,3,4\}$,
we have
$\beta_{1}^{+}+\beta_{2}^{+}+\beta_{3}^{+}+\beta_{4}^{+}\in (2,3)$.
It implies that $q\in(1,\infty)$.
By Young's inequality,
\begin{align}\label{result:case2}
\left\|\int_{\mathbb{R}}
  \frac{ I^{1/2}_{\beta_{4}}|\widehat{\psi}|^{2}(u)}
{|u|^{1-\beta_{1}}}
\frac{Q(\bullet-u)}{|\bullet-u|^{1-\beta_{3}}}
\,du\right\|_{q}\leq
 \left\|\frac{I^{1/2}_{\beta_{4}}|\widehat{\psi}|^{2}(\bullet)}{|\bullet|^{1-\beta_{1}}}\right\|_{p}
 \left\|\frac{Q(\bullet)}{|\bullet|^{1-\beta_{3}}}\right\|_{r}.
\end{align}

The assertion $\mathcal{K}<\infty$
follows from (\ref{result:case1}), (\ref{result:case2}), the assumption $\psi, \widehat{\psi}\in L^{1}$, and H$\ddot{\textup{o}}$lder's inequality. In the following, we prove \ref{Claim 1}, \ref{Claim 2}, and \ref{Claim 3}.

{\it Proof of \ref{Claim 1}:} First of all, according to the Hardy-Littlewood-Sobolev theorem for fractional integration
\cite[p. 119]{stein1970singular}, for any
$1<p'<q'<\infty$ and $\beta\in(0,1)$ with $1/q' = 1/p'-\beta$, there exists a constant $A_{p',q'}$ such that
\begin{align*}
\| I_{\beta}|\widehat{\psi}|^{2}\|_{q'}\leq A_{p',q'}
\||\widehat{\psi}|^{2}\|_{p'}.
\end{align*}
Because $\||\widehat{\psi}|^{2}\|_{p'}<\infty$ for any $p'>1$,
we have $I_{\beta}|\widehat{\psi}|^{2}\in L^{1/(1-\beta^{+})}$, where $\beta^{+}$ is a constant greater than $\beta$, sufficiently proximate to it.
Additionally, for any $w,w'\in \mathbb{R}$, $\beta\in(1/2,1)$, and $\widehat{\psi}\in L^{1}\cap L^{\infty}$,
\begin{align*}
\left|I_{\beta}|\widehat{\psi}|^{2}(w)-I_{\beta}|\widehat{\psi}|^{2}(w')\right|
=&
\left|\int_{\mathbb{R}}
 \frac{1}{
|z|^{1-\beta}}\left[|\widehat{\psi}(w-z)|^{2}-|\widehat{\psi}(w'-z)|^{2}\right]
dz\right|
\\\leq&
\int_{-1}^{1}
\frac{1}{
|z|^{1-\beta}}\left||\widehat{\psi}(w-z)|^{2}-|\widehat{\psi}(w'-z)|^{2}\right|
dz
\\&+
\int_{\mathbb{R}}
1_{\{|z|>1\}}\frac{1}{
|z|^{1-\beta}}\left||\widehat{\psi}(w-z)|^{2}-|\widehat{\psi}(w'-z)|^{2}\right|
dz
\\\leq&
\left[\int_{-1}^{1}
\frac{dz}{
|z|^{2(1-\beta)}}\right]^{1/2}
\left[\int_{-1}^{1}
\left||\widehat{\psi}(w-z)|^{2}-|\widehat{\psi}(w'-z)|^{2}\right|^{2}
dz\right]^{1/2}
\\&+
\int_{\mathbb{R}}
\left||\widehat{\psi}(w-z)|^{2}-|\widehat{\psi}(w'-z)|^{2}\right|
dz,
\end{align*}
which tends to zero as $|w-w'|\rightarrow 0$ due to the continuity of translation in the $L^{r'}$ norm, where $r'\in[1,\infty)$.
Hence, for any $\beta\in(1/2,1)$,
\begin{equation}\label{I_beta_L_C}
I^{1/2}_{\beta}|\widehat{\psi}|^{2}\in L^{2/(1-\beta^{+})} \cap C(\mathbb{R}).
\end{equation}
Because $I^{1/2}_{\beta_{4}}|\widehat{\psi}|^{2}(\cdot)$ is an even function,
\begin{align*}
\int_{\mathbb{R}}
\left(\frac{I^{1/2}_{\beta_{4}}|\widehat{\psi}|^{2}(w)}{|w|^{1-\beta_{1}}}\right)^{2/(3-2\beta_{1}^{+}-\beta_{4}^{+})}dw
=&
2\int_{0}^{1}
\left(\frac{I^{1/2}_{\beta_{4}}|\widehat{\psi}|^{2}(w)}{|w|^{1-\beta_{1}}}\right)^{2/(3-2\beta_{1}^{+}-\beta_{4}^{+})}dw
\\+&
2\int_{1}^{\infty}
\left(\frac{I^{1/2}_{\beta_{4}}|\widehat{\psi}|^{2}(w)}{|w|^{1-\beta_{1}}}\right)^{2/(3-2\beta_{1}^{+}-\beta_{4}^{+})}dw
:=\mathcal{I}_{1}+\mathcal{I}_{2}.
\end{align*}
Because $2(1-\beta_{1})/(3-2\beta_{1}^{+}-\beta_{4}^{+})\in (0,1)$ and (\ref{I_beta_L_C}) implies the boundedness of $I^{1/2}_{\beta_{4}}|\widehat{\psi}|^{2}$, we have
\begin{align*}
\mathcal{I}_{1} \lesssim
\int_{0}^{1}
|w|^{-\frac{2(1-\beta_{1})}{3-2\beta_{1}^{+}-\beta_{4}^{+}}}dw<\infty.
\end{align*}
Regarding the estimation of $\mathcal{I}_{2}$, because
\begin{align*}
\left(\frac{3-2\beta_{1}^{+}-\beta_{4}^{+}}{1-\beta_{4}^{+}}\right)^{-1}+
\left(\frac{3-2\beta_{1}^{+}-\beta_{4}^{+}}{2(1-\beta_{1}^{+})}\right)^{-1}
=1
\end{align*}
and
$$
\frac{3-2\beta_{1}^{+}-\beta_{4}^{+}}{1-\beta_{4}^{+}}>\frac{3-2\beta_{1}^{+}-\beta_{4}^{+}}{2(1-\beta_{1}^{+})}>1
$$
for $\beta\in (1/2,3/4)$,
we apply H$\ddot{\textup{o}}$lder's inequality to derive
\begin{align*}
\mathcal{I}_{2}=&2\int_{1}^{\infty}
\left(\frac{I^{1/2}_{\beta_{4}}|\widehat{\psi}|^{2}(w)}{|w|^{1-\beta_{1}}}\right)^{2/(3-2\beta_{1}^{+}-\beta_{4}^{+})}dw
\\\lesssim&
\left[\int_{1}^{\infty}
\left|I_{\beta_{4}}|\widehat{\psi}|^{2}(w)\right|^{1/(1-\beta_{4}^{+})}dw\right]^{(1-\beta_{4}^{+})/(3-2\beta_{1}^{+}-\beta_{4}^{+})}
\left[\int_{1}^{\infty}
|w|^{-\frac{1-\beta_{1}}{1-\beta_{1}^{+}}}dw
\right]^{2(1-\beta_{1}^{+})/(3-2\beta_{1}^{+}-\beta_{4}^{+})}.
\end{align*}
The convergence of the right-hand side of this inequality stems from (\ref{I_beta_L_C}).
Consequently, the proof of  \ref{Claim 1} is concluded.

{\it Proof of \ref{Claim 2}:} Referring to the definition of $Q$ in (\ref{def:Q}) and utilizing Minkowski's inequality for integrals, we derive
\begin{align}\notag
\|Q\|_{2/(3-2\beta_{2}^{+}-\beta_{4}^{+})}\leq
\|\widehat{\psi}\|_{1}
\left\| \frac{I^{1/2}_{\beta_{4}}|\widehat{\psi}|^{2}(\bullet)}{
|\bullet|^{1-\beta_{2}}}\right\|_{2/(3-2\beta_{2}^{+}-\beta_{4}^{+})}<\infty,
\end{align}
where the final inequality stems from \ref{Claim 1} by substituting $\beta_{1}$ with $\beta_{2}$.
Because $\beta\in(1/2,3/4)$, both
$2/(3-2\beta_{2}^{+}-\beta_{4}^{+})$ and $2/(2\beta_{2}^{+}+\beta_{4}^{+}-1)$ are strictly greater than one.
Moreover, the following equality
$$
\left(\frac{2}{3-2\beta_{2}^{+}-\beta_{4}^{+}}\right)^{-1}+\left(\frac{2}{2\beta_{2}^{+}+\beta_{4}^{+}-1}\right)^{-1}=1
$$
holds. It allows us to apply H$\ddot{\textup{o}}$lder's inequality. For any $w,w'\in \mathbb{R}$,
\begin{align}\notag
\left|Q(w)-Q(w')\right|\leq& \int_{\mathbb{R}}\left|\widehat{\psi}(y)\right|
\left| \frac{I^{1/2}_{\beta_{4}}|\widehat{\psi}|^{2}(w+y)}{
|w+y|^{1-\beta_{2}}}-\frac{I^{1/2}_{\beta_{4}}|\widehat{\psi}|^{2}(w'+y)}{
|w'+y|^{1-\beta_{2}}}\right|
\,dy
\\\notag\leq& \|\widehat{\psi}\|_{2/(2\beta_{2}^{+}+\beta_{4}^{+}-1)}
\left\| \frac{I^{1/2}_{\beta_{4}}|\widehat{\psi}|^{2}(w+\bullet)}{
|w+\bullet|^{1-\beta_{2}}}-\frac{I^{1/2}_{\beta_{4}}|\widehat{\psi}|^{2}(w'+\bullet)}{
|w'+\bullet|^{1-\beta_{2}}}\right\|_{2/(3-2\beta_{2}^{+}-\beta_{4}^{+})},
\end{align}
which tends to zero as $|w-w'|\rightarrow 0$ due to the continuity of translation in the $L^{r'}$ norm, where $r'\in[1,\infty)$.
Hence, $Q$ is continuous everywhere, thereby concluding the proof of \ref{Claim 2}.

{\it Proof of \ref{Claim 3}:}
Because $Q$ is continuous everywhere,
\begin{equation*}
\int_{\mathbb{R}}\left|\frac{Q(w)}{|w|^{1-\beta_{3}}}\right|^{2/(5-2\beta_{2}^{+}-2\beta_{3}^{+}-\beta_{4}^{+})}dw
\lesssim \mathcal{I}_{3}+\mathcal{I}_{4},
\end{equation*}
where
\begin{equation*}
\mathcal{I}_{3} = \int_{0}^{1}
|w|^{-\frac{2(1-\beta_{3})}{5-2\beta_{2}^{+}-2\beta_{3}^{+}-\beta_{4}^{+}}}dw
\end{equation*}
and
\begin{equation*}
\mathcal{I}_{4} = \int_{1}^{\infty}
\left|\frac{Q(w)}{|w|^{1-\beta_{3}}}\right|^{2/(5-2\beta_{2}^{+}-2\beta_{3}^{+}-\beta_{4}^{+})}dw.
\end{equation*}
It is clear that $\mathcal{I}_{3}<\infty$.
Because
\begin{equation*}
\left(\frac{5-2\beta_{2}^{+}-2\beta_{3}^{+}-\beta_{4}^{+}}{3-2\beta_{2}^{+}-\beta_{4}^{+}}\right)^{-1}
+\left(\frac{5-2\beta_{2}^{+}-2\beta_{3}^{+}-\beta_{4}^{+}}{2(1-\beta_{3}^{+})}\right)^{-1}
=1
\end{equation*}
and both terms within the parentheses are greater than one,
we have
\begin{align}\label{final:proof:claim3}
\mathcal{I}_{4} \leq \left[\int_{1}^{\infty}
|Q(w)|^{2/(3-2\beta_{2}^{+}-\beta_{4}^{+})}dw\right]^{\frac{3-2\beta_{2}^{+}-\beta_{4}^{+}}{5-2\beta_{2}^{+}-2\beta_{3}^{+}-\beta_{4}^{+}}}
\left[\int_{1}^{\infty}
|w|^{\frac{1-\beta_{3}}{1-\beta_{3}^{+}}}dw \right]^{\frac{2(1-\beta_{3}^{+})}{5-2\beta_{2}^{+}-2\beta_{3}^{+}-\beta_{4}^{+}}}.
\end{align}
The convergence of the first integral in (\ref{final:proof:claim3}) follows from \ref{Claim 2},  thereby establishing the proof of \ref{Claim 3}.

\subsection{Proof of Corollary \ref{corollary:cov:chi}}\label{sec:proof:corollary:cov:chi}
Denote the limiting process of $\{sU^{2}[s]Y(st)\}_{t\in \mathbb{R}}$ by $\{F(t)\}_{t\in \mathbb{R}}$.
Consider two mean-zero Gaussian processes, denoted by $\{N_{R}(t)\}_{t\in \mathbb{R}}$ and $\{N_{I}(t)\}_{t\in \mathbb{R}}$, with the following covariance matrix:
\begin{align}\notag
&\mathbb{E}\left[\begin{array}{cc}
N_{R}(t)N_{R}(t')  & N_{R}(t)N_{I}(t')\\
N_{I}(t)N_{R}(t') & N_{I}(t)N_{I}(t')
\end{array}\right]
\\\label{def:cov:NRNI}=&
2^{4}\|f_{X}\|_{2}^{2}\int_{0}^{\infty}
\hspace{-0.1cm}\left[\hspace{-0.1cm}\begin{array}{cc}
\cos((t-t')u)   & \sin((t-t')u)
\\
\sin((t'-t)u) &  \cos((t-t')u)
\end{array}\hspace{-0.1cm}\right]
|\widehat{\psi}(u)|^{2}du,
\end{align}
where $t,t'\in \mathbb{R}.$
By Theorem \ref{thm:Stein},
$F(t) \overset{d}{=} N^{2}_{R}(t)+N^{2}_{I}(t),$
where
$\overset{d}{=}$ means that the finite-dimensional distributions of $F$ and $N^{2}_{R}+N^{2}_{I}$ are identical.
From (\ref{def:cov:NRNI}), we obtain that
\begin{align}\notag
&\textup{Cov}\left(F(t),F(t')\right) = \textup{Cov}\left(N^{2}_{R}(t)+N^{2}_{I}(t),N^{2}_{R}(t')+N^{2}_{I}(t')\right)
\\\notag=&\textup{Cov}\left(N^{2}_{R}(t),N^{2}_{R}(t')\right)+
\textup{Cov}\left(N^{2}_{R}(t),N^{2}_{I}(t')\right)
+\textup{Cov}\left(N^{2}_{I}(t),N^{2}_{R}(t')\right)
+
 \textup{Cov}\left(N^{2}_{I}(t),N^{2}_{I}(t')\right)
\\\notag=&
2\left[\textup{Cov}\left(N_{R}(t),N_{R}(t')\right)^{2}
+\textup{Cov}\left(N_{R}(t),N_{I}(t')\right)^{2}
+\textup{Cov}\left(N_{I}(t),N_{R}(t')\right)^{2}
+
\textup{Cov}\left(N_{I}(t),N_{I}(t')\right)^{2}\right]
\\\notag=&
2^{10}\|f_{X}\|_{2}^{4}\left(\int_{0}^{\infty}\cos((t-t')u)|\widehat{\psi}(u)|^{2}du\right)^{2}
+2^{10}\|f_{X}\|_{2}^{4}\left(\int_{0}^{\infty}\sin((t-t')u)|\widehat{\psi}(u)|^{2}du\right)^{2}.
\end{align}



\section*{Acknowledgments}
The author expresses gratitude to the editors for handling this paper and to the anonymous reviewers for their valuable comments,
which have significantly improved its quality. The author also thanks Hau-Tieng Wu and Yuan-Chung Sheu for their guidance. Additionally, appreciation is extended for the financial support provided by the National Center for Theoretical Sciences
and the National Science and Technology Council, Taiwan (Project No. 110-2628-M-006-003-MY3).




\appendix

\section{Lemmas about the slowly varying function $L$}\label{sec:appendix:L}
In the proof of Theorem \ref{thm:nonStein}, the estimate of
$$
\left|\frac{L^{2}\left(s/|x|\right)}{L^{2}(s)}-1\right|
$$
for the two cases $|x|\leq1$ and $|x|>1$ is required. Below, we integrate them for ease of reference and application.

\begin{Lemma}\label{lemma:L:ratio:easy}
Under Assumption \ref{Assumption:spectral},
for any $\delta>0$ and sufficiently large $s$, there exists a constant $C>0$ such that the following inequality
\begin{align*}
\left|\frac{L^{2}\left(s/|x|\right)}{L^{2}(s)}-1\right|\leq C
g_{\tau}(s)|x|^{-\delta}
\end{align*}
holds for any $x\in \mathbb{R}\setminus\{0\}$ with $|x|\leq 1$.
\end{Lemma}

{\it Proof:}
Under Assumption \ref{Assumption:spectral},
there exists a constant $C_{1}$ such that for any $s>0$ and $|x|\leq1$,
\begin{align}\notag
\left|\frac{L^{2}\left(s/|x|\right)}{L^{2}(s)}-1\right|
=&
\left|\frac{L\left(s/|x|\right)}{L(s)}-1\right|\left|\frac{L\left(s/|x|\right)}{L(s)}+1\right|
\\\notag\leq& C_{1} g_{\tau}(s)h_{\tau}\left(\frac{1}{|x|}\right)\left[C_{1} g_{\tau}(s)h_{\tau}\left(\frac{1}{|x|}\right)+2\right]
\\\label{easy1}=&2C_{1} g_{\tau}(s)h_{\tau}\left(\frac{1}{|x|}\right)+C_{1}^{2}g^{2}_{\tau}(s)h^{2}_{\tau}\left(\frac{1}{|x|}\right).
\end{align}
We estimate the right hand side of (\ref{easy1}) by considering two cases as follows.
\begin{itemize}
\item For $\tau=0$, $h_{\tau}\left(|x|^{-1}\right) = \ln(|x|^{-1})$.
For any $\delta>0$, there exists a constant $C_{2}>0$ such that
$\ln\left(|x|^{-1}\right)\leq C_{2}|x|^{-\delta}$ for all $|x|\in (0,1]$. Consequently, (\ref{easy1}) implies that
\begin{align}\label{easy2.1}
\left|\frac{L^{2}\left(s/|x|\right)}{L^{2}(s)}-1\right|
\leq 2C_{1}C_{2} g_{\tau}(s)|x|^{-\delta}+C_{1}^{2}C_{2}^{2}g^{2}_{\tau}(s)|x|^{-2\delta}.
\end{align}

\item For $\tau<0$, the definition of $h_{\tau}$ in (\ref{def:h_tau}) implies that $h_{\tau}(\cdot)\leq |\tau|^{-1}$.
Hence,
\begin{align}\label{easy2.2}
\left|\frac{L^{2}\left(s/|x|\right)}{L^{2}(s)}-1\right|
\leq 2C_{1}|\tau|^{-1} g_{\tau}(s)+C_{1}^{2}|\tau|^{-2}g^{2}_{\tau}(s).
\end{align}
\end{itemize}
The proof is completed by individually applying the assumption $g_{\tau}(s)\rightarrow 0$ as $s\rightarrow\infty$ to (\ref{easy2.1}) and (\ref{easy2.2}).
\qed

\begin{Lemma}\label{lemma:L:ratio:diff}
Under Assumption \ref{Assumption:spectral},
for any $\delta>0$ and sufficiently large $s$, there exists a constant $C>0$ such that the following inequality
\begin{align}\label{lemma:L:ratio:diff:eq}
\left|\frac{L^{2}\left(s/|x|\right)}{L^{2}(s)}-1\right|\leq C g_{\tau}(s)|x|^{-2\tau+\delta}
\end{align}
holds for any $x\in \mathbb{R}$ with $|x|>1$.
\end{Lemma}

{\it Proof:}
First of all, we rewrite the left hand side of (\ref{lemma:L:ratio:diff:eq}) as follows
\begin{align}\label{proof:slowly:diff1}
\left|\frac{L^{2}\left(s/|x|\right)}{L^{2}(s)}-1\right|
=\frac{L^{2}\left(s/|x|\right)}{L^{2}(s)}\left|1-\frac{L^{2}(s)}{L^{2}\left(s/|x|\right)}\right|.
\end{align}
From \cite[Theorem 1.5.3]{bingham1989regular}, we know that if a mapping $s\mapsto R_{\rho}(s)$
varies regularly with index $\rho>0$ and
it is locally bounded on $[0,\infty)$,
we have
\begin{align}\label{cited_regular_theorem}
\underset{s\rightarrow\infty}{\lim}\ \underset{0\leq r\leq s }{\sup}\ \frac{R_{\rho}\left(r\right)}{R_{\rho}(s)}=1.
\end{align}
For any $\delta>0$, by applying (\ref{cited_regular_theorem}) to the function
$s\mapsto  s^{\delta}L^{2}(s)$, which varies regularly with index $\delta$,
there exist a constant $C_{1}>0$ and a sufficiently large $s_{0}$ such that
\begin{align}\label{proof:slowly:diff3}
\frac{L^{2}\left(s/|x|\right)}{L^{2}(s)}
\leq&
\left[\underset{1<|z|<\infty }{\sup} \frac{(s/|z|)^{\delta}L^{2}\left(s/|z|\right)}{s^{\delta}L^{2}(s)}\right]|x|^{\delta}
\leq C_{1}|x|^{\delta}
\end{align}
for any $x\in \mathbb{R}$ with $|x|>1$ and $s\geq s_{0}$.
For the second term in (\ref{proof:slowly:diff1}), by Assumption \ref{Assumption:spectral},
there exists a constant $C_{2}>0$ such that
\begin{align}\notag
\left|1-\frac{L^{2}(s)}{L^{2}\left(s/|x|\right)}\right|
=& \left|1-\frac{L(s)}{L\left(s/|x|\right)}\right|\left|1+\frac{L(s)}{L\left(s/|x|\right)}\right|
\\\notag\leq&  C_{2}g_{\tau}(s/|x|)h_{\tau}(|x|)\left[2+C_{2}g_{\tau}(s/|x|)h_{\tau}(|x|)\right]
\\\label{challenge1}=& 2C_{2}g_{\tau}(s/|x|)h_{\tau}(|x|)+C_{2}^{2}g_{\tau}^{2}(s/|x|)h^{2}_{\tau}(|x|)
\end{align}
for all $s>0$ and $|x|>1$.
For the term $g_{\tau}(s/|x|)$ in (\ref{challenge1}), we note that for $|x|>1$,
\begin{align*}
g_{\tau}(s/|x|) =& g_{\tau}(s)\frac{g_{\tau}(s/|x|)}{g_{\tau}(s)}
\\\leq&
g_{\tau}(s) \left[\underset{1<|z|<\infty }{\sup}\frac{(s/|z|)^{-\tau+\delta}g_{\tau}(s/|z|)}{s^{-\tau+\delta}g_{\tau}(s)}\right]|x|^{-\tau+\delta}.
\end{align*}
By applying (\ref{cited_regular_theorem}) to the function
$s\mapsto  s^{-\tau+\delta}g_{\tau}(s)$ inside the bracket above, which varies regularly with index $\delta$,
we obtain that there exist a constant $C_{3}$ and a sufficiently large $s_{0}$ such that
\begin{align}\label{estimate:g(s/x)}
g_{\tau}(s/|x|)
\leq C_{3}g_{\tau}(s) |x|^{-\tau+\delta}
\end{align}
holds for $s\geq s_{0}$ and $|x|>1$.
Substituting (\ref{estimate:g(s/x)}) into (\ref{challenge1}) yields
\begin{align}\label{proof:slowly:diff2}
\left|1-\frac{L^{2}(s)}{L^{2}\left(s/|x|\right)}\right|
\leq 2C_{2}C_{3}g_{\tau}(s) |x|^{-\tau+\delta}h_{\tau}(|x|)+C_{2}^{2}C_{3}^{2}g_{\tau}^{2}(s) |x|^{-2\tau+2\delta}h^{2}_{\tau}(|x|)
\end{align}
for $s\geq s_{0}$ and $|x|>1$.
By substituting (\ref{proof:slowly:diff3}) and (\ref{proof:slowly:diff2}) into
(\ref{proof:slowly:diff1}), for sufficiently large $s$, there exists a constant $C_{4}$ such that
\begin{align}\label{proof:slowly:diff4}
\left|\frac{L^{2}\left(s/|x|\right)}{L^{2}(s)}-1\right|
\leq C_{4}\left[ g_{\tau}(s) |x|^{-\tau+2\delta}h_{\tau}(|x|)+g_{\tau}^{2}(s) |x|^{-2\tau+3\delta}h^{2}_{\tau}(|x|)\right]
\end{align}
for $s\geq s_{0}$ and $|x|>1$. Finally, we consider two cases based on the definition of $h_{\tau}(\cdot)$ in (\ref{def:h_tau}).
\begin{itemize}
\item For $\tau=0$, $h_{\tau}\left(|x|\right) = \ln(|x|)$.
There exists a constant $C_{5}$ such that
$\ln(|x|)\leq C_{5}|x|^{\delta}$ for all $|x|>1$. Consequently, (\ref{proof:slowly:diff4}) implies that
there exists a constant $C_{6}$ such that
\begin{align}\label{diff2.1}
\left|\frac{L^{2}\left(s/|x|\right)}{L^{2}(s)}-1\right|
\leq C_{6}\left[ g_{\tau}(s) |x|^{-\tau+3\delta}+g_{\tau}^{2}(s) |x|^{-2\tau+5\delta}\right]
\end{align}
for any $s\geq s_{0}$ and $|x|>1$.

\item For $\tau<0$, $h_{\tau}(\cdot)\leq |\tau|^{-1}$.
Hence, (\ref{proof:slowly:diff4}) implies that
\begin{align}\label{diff2.2}
\left|\frac{L^{2}\left(s/|x|\right)}{L^{2}(s)}-1\right|
\leq C_{4}\left[|\tau|^{-1} g_{\tau}(s) |x|^{-\tau+2\delta}+|\tau|^{-2}g_{\tau}^{2}(s) |x|^{-2\tau+3\delta}\right].
\end{align}
\end{itemize}
The proof is completed by individually applying the assumption $g_{\tau}(s)\rightarrow 0$ as $s\rightarrow\infty$ to (\ref{diff2.1}) and (\ref{diff2.2}).
\qed

\section{Proof of the anti-concentration inequality (\ref{levyineq})}\label{sec:proof:levyineq}
The following proof follows the methodology used in the proofs of Theorem 3.1 and Lemma 4.3 in \cite{MR3003367},
as well as the proof of Lemma 4 in \cite{MR4003569}.

Let us denote by $I_{k}(\cdot)$ a multiple Wiener-It$\hat{\textup{o}}$ integral of order $k$, where $k\in \mathbb{N}$.
From the definition of $Z$ in (\ref{def:Z0})-(\ref{def:Zlimit2}),
\begin{align}
Z(t) = Z^{(0)}+I_{2}(g^{(2)})+I_{4}(g^{(4)}),
\end{align}
where $Z^{(0)}$ is a non-random constant,
\begin{align*}
g^{(2)}(\lambda_{1},\lambda_{2}) =&
8
e^{it(\lambda_{1}+\lambda_{2})}\left|\lambda_{1}\lambda_{2}\right|^{\frac{\beta-1}{2}}
1_{\{\lambda_{1}\lambda_{2}<0\}}
\\&\times\int_{\mathbb{R}}
\frac{\widehat{\psi}(\lambda_{1}+\xi)\widehat{\psi}(\lambda_{2}-\xi)}{\left|\xi\right|^{1-\beta}} 1_{\{(\lambda_{1}+\xi)(\lambda_{2}-\xi)<0\}}d\xi
\end{align*}
and $\widetilde{g}^{(4)}$ is the symmetrization of  $g^{(4)}$, while $g^{(4)}$ is defined as
\begin{align*}
g^{(4)}(\lambda_{1},\lambda_{2},\lambda_{3},\lambda_{4}) =&8e^{it(\lambda_{1}+\lambda_{2}+\lambda_{3}+\lambda_{4})}\frac{\widehat{\psi}(\lambda_{1}+\lambda_{2})
\widehat{\psi}(\lambda_{3}+\lambda_{4})}{|\lambda_{1}\lambda_{2}\lambda_{3}\lambda_{4}|^{(1-\beta)/2}}
\\&\times1_{\{\lambda_{1}\lambda_{2}<0\}}1_{\{\lambda_{3}\lambda_{4}<0\}}
1_{\{(\lambda_{1}+\lambda_{2})(\lambda_{3}+\lambda_{4})<0\}}.
\end{align*}
By Lemma \ref{lemma:integral:nu12}, $g^{(2)}\in L^{2}(\mathbb{R}^{2})$. On the other hand, because $\beta\in (0,1/2)$ under the hypothesis of Theorem \ref{thm:nonStein},
the Riesz composition formula shows that $g^{(4)}\in L^{2}(\mathbb{R}^{4})$ with $\|\widetilde{g}^{(4)}\|_{2}\leq \|g^{(4)}\|_{2}.$
Let $\{e_{i}\}_{i\in \mathbb{N}}$ be an orthonormal basis of $L^{2}(\mathbb{R})$.
Then, there exist constants $\{c_{i_{1},i_{2}}\}_{i_{1},i_{2}\in \mathbb{N}}$ and $\{d_{i_{1},i_{2},i_{3},i_{4}}\}_{i_{1},i_{2},i_{3},i_{4}\in \mathbb{N}}$ such that
\begin{align*}
g^{(2)} =&
\underset{i_{1},i_{2}\in \mathbb{N}}{\sum}c_{i_{1},i_{2}} e_{i_{1}}\otimes e_{i_{2}}
\end{align*}
and
\begin{align*}
\widetilde{g}^{(4)} =&
\underset{i_{1},i_{2},i_{3},i_{4}\in \mathbb{N}}{\sum}d_{i_{1},i_{2},i_{3},i_{4}} e_{i_{1}}\otimes e_{i_{2}}\otimes e_{i_{3}}\otimes e_{i_{4}}.
\end{align*}
For any $n\in \mathbb{N}$, denote
\begin{align*}
g^{(2)}_{n} =&
\underset{i_{1},i_{2}\in \{1,2,...,n\}}{\sum}c_{i_{1},i_{2}} e_{i_{1}}\otimes e_{i_{2}},
\end{align*}
\begin{align*}
\widetilde{g}_{n}^{(4)} =&
\underset{i_{1},i_{2},i_{3},i_{4}\in\{1,2,...,n\}}{\sum}d_{i_{1},i_{2},i_{3},i_{4}} e_{i_{1}}\otimes e_{i_{2}}\otimes e_{i_{3}}\otimes e_{i_{4}},
\end{align*}
and
\begin{align*}
Z(t;n)=Z^{(0)}+I_{2}(g^{(2)}_{n})+I_{4}(\widetilde{g}_{n}^{(4)}).
\end{align*}
Because $\|g^{(2)}_{n}-g^{(2)}\|_{2}$ and $\|\widetilde{g}_{n}^{(4)}-\widetilde{g}^{(4)}\|_{2}$ converge to zero
when $n\rightarrow\infty$, the isometry property of Wiener-It$\hat{\textup{o}}$ integrals implies that
$\mathbb{E}|Z(t;n)- Z(t)|^{2}\rightarrow0$ when $n\rightarrow\infty$.
Hence, there exists a subsequence $\{n_{j}\}_{j\in \mathbb{N}}$ of $\{1,2,3,...\}$ such that
$Z(t;n_{j})\rightarrow Z(t)$ almost surely
when $j\rightarrow\infty$.
It also follows that for any $\varepsilon>0$ and $z\in \mathbb{R}$,
\begin{align}\notag
&\mathbb{P}\left(|Z(t)-z|\leq \varepsilon/2\right) =  \mathbb{P}\left(\underset{j\rightarrow\infty}{\liminf}\  |Z(t;n_{j})-z|\leq \varepsilon/2\right)
\\\notag\leq& \mathbb{E}\left[\underset{j\rightarrow\infty}{\liminf}\ 1_{[0,\varepsilon]}\left(|Z(t;n_{j})-z|\right)\right]
\leq \underset{j\rightarrow\infty}{\liminf}\ \mathbb{E}\left[1_{[0,\varepsilon]}\left(|Z(t;n_{j})-z|\right)\right]
\\\label{appendix:fatou}=&\underset{j\rightarrow\infty}{\liminf}\ \mathbb{P}\left(|Z(t;n_{j})-z|\leq \varepsilon\right).
\end{align}

For estimating the right hand side of (\ref{appendix:fatou}),
we note that
\begin{align}\notag
I_{2}(g^{(2)}_{n_{j}}) =& I_{2}\left(\underset{i_{1},i_{2}\in \{1,2,...,n_{j}\}}{\sum}c_{i_{1},i_{2}} e_{i_{1}}\otimes e_{i_{2}}
\right)
\\\notag=& \underset{i_{1},i_{2}\in \{1,2,...,n_{j}\}}{\sum}c_{i_{1},i_{2}} I_{2}\left(e_{i_{1}}\otimes  e_{i_{2}}
\right),
\end{align}
where, by It$\hat{\textup{o}}$'s formula \cite[Theorem 4.3]{major1981lecture},
\begin{align*}
I_{2}\left(e_{i_{1}}\otimes  e_{i_{2}}
\right) =\left\{\begin{array}{ll} H_{1}\left(\int_{\mathbb{R}}e_{i_{1}}(\lambda)W(d\lambda)\right) H_{1}\left(\int_{\mathbb{R}}e_{i_{2}}(\lambda)W(d\lambda)\right) & \textup{if}\ i_{1}\neq i_{2},
\\
H_{2}\left(\int_{\mathbb{R}}e_{i_{1}}(\lambda)W(d\lambda)\right)& \textup{if}\ i_{1}= i_{2}.
\end{array}\right.
\end{align*}
Hence,
\begin{align}\label{appendix:poly:I2}
I_{2}(g^{(2)}_{n_{j}}) = U_{2,n_{j}}\left(\int_{\mathbb{R}}e_{1}(\lambda)W(d\lambda),...,\int_{\mathbb{R}}e_{n_{j}}(\lambda)W(d\lambda)\right),
\end{align}
for a certain polynomial function $U_{2,n_{j}}: \mathbb{R}^{n_{j}}\rightarrow \mathbb{R}$ of degree two.
Similarly,
\begin{align}\label{appendix:poly:I4}
I_{4}(\widetilde{g}^{(4)}_{n_{j}}) = U_{4,n_{j}}\left(\int_{\mathbb{R}}e_{1}(\lambda)W(d\lambda),...,\int_{\mathbb{R}}e_{n_{j}}(\lambda)W(d\lambda)\right)
\end{align}
for a certain polynomial function $U_{4,n_{j}}: \mathbb{R}^{n_{j}}\rightarrow \mathbb{R}$ of degree four.
For any fixed $z\in \mathbb{R}$, let $V_{4,n_{j}}= -z+Z^{(0)}+U_{2,n_{j}}+U_{4,n_{j}}$, which is a polynomial function of degree four.
By (\ref{appendix:poly:I2}) and (\ref{appendix:poly:I4}),
\begin{align}\label{appendix:poly:Zn}
Z(t;n_{j})-z=V_{4,n_{j}}\left(\int_{\mathbb{R}}e_{1}(\lambda)W(d\lambda),...,\int_{\mathbb{R}}e_{n_{j}}(\lambda)W(d\lambda)\right).
\end{align}
Because the inputs of the polynomial function $V_{4,n_{j}}$ in (\ref{appendix:poly:Zn}) consist of $n_{j}$
independent standard normal random variables, by the Carbery-Wright inequality \cite[Theorem 2.5]{MR3003367},
there exists a constant $C$ independent of $n_{j}$ such that
\begin{align}\label{appendix:Carbery}
\mathbb{P}\left(|Z(t;n_{j})-z|\leq \varepsilon\right)\leq C \left\{\mathbb{E}\left[\left|Z(t;n_{j})-z\right|^{2}\right]\right\}^{-\frac{1}{8}}\varepsilon^{\frac{1}{4}}.
\end{align}
Since $\mathbb{E}|Z(t;n_{j})- Z(t)|^{2}\rightarrow0$ when $j\rightarrow\infty$,
\begin{align}\label{appendix:fatou2}
\underset{j\rightarrow\infty}{\liminf}\left\{\mathbb{E}\left[\left|Z(t;n_{j})-z\right|^{2}\right]\right\}^{-\frac{1}{8}}
=\left\{\mathbb{E}\left[\left|Z(t)-z\right|^{2}\right]\right\}^{-\frac{1}{8}}.
\end{align}
Combining (\ref{appendix:Carbery}) and (\ref{appendix:fatou2}) yields
\begin{align}\label{appendix:fatou3}
\underset{j\rightarrow\infty}{\liminf}\ \mathbb{P}\left(|Z(t;n_{j})-z|\leq \varepsilon\right)
\leq C\left\{\mathbb{E}\left[\left|Z(t)-z\right|^{2}\right]\right\}^{-\frac{1}{8}}\varepsilon^{\frac{1}{4}}.
\end{align}
By substituting (\ref{appendix:fatou3}) into  (\ref{appendix:fatou}), we obtain the desired inequality
\begin{align}\notag
\mathbb{P}\left(|Z(t)-z|\leq \varepsilon/2\right) \leq C\left\{\mathbb{E}\left[\left|Z(t)-z\right|^{2}\right]\right\}^{-\frac{1}{8}}\varepsilon^{\frac{1}{4}}.
\end{align}
\qed

\end{document}